\documentclass[11pt]{ip-journal}
\usepackage[includefoot]{geometry}
\usepackage{color}
\usepackage{amssymb}
\usepackage{amsmath}
\usepackage{amsthm} 
\usepackage{mathrsfs}                                                            
\usepackage[all]{xy}
\usepackage{graphicx}
\usepackage{wrapfig}

\def\2{{1\over 2}}

\newcommand{\rf}[1]{(\ref{#1})}
\def\b{\bar}

\renewcommand{\t}{\tilde}
\newcommand{\p}{\partial}

\def\p{\partial}

\def\c1m2{\sqrt{{\bf x_1}}-c \sqrt{\bf x_2}}

\def\b{\bar}

\def\<{\langle}
\def\>{\rangle}
\def\+{\dagger}

\newtheorem{Thm}{Theorem}[section]
\newtheorem{Lem}[Thm]{Lemma}
\newtheorem{Prop}[Thm]{Proposition}
\newtheorem{Cor}[Thm]{Corollary}

\begin{document}
\title{Decorated Super-Teichm\"uller Space}
\author{R.\ C.\ Penner}
\address{
R.\ C.\ Penner,\newline
Institut des Hautes \'Etudes Scientifiques\newline
Le Bois-Marie 35, route de Chartres\newline
91440 Bures-sur-Yvette, France\newline
rpenner@ihes.fr
\newline\newline
Department of Mathematics \newline
University of California--Los Angeles\newline
Box 951555\newline
Los Angeles, CA 90095 USA}
\author{Anton M. Zeitlin}
\address{ Anton M. Zeitlin,\newline
Department of Mathematics,\newline
303 Lockett Hall\newline
Louisiana State University,\newline
Baton Rouge, LA 70803 USA
\newline\newline
IPME RAS, 
V.O. Bolshoj pr., 61,\newline 
199178, St. Petersburg\newline
zeitlin@lsu.edu,\newline
http://math.lsu.edu/$\sim$zeitlin \newline
http://www.ipme.ru/zam.html  }

\begin{abstract}  We introduce coordinates for a principal bundle $S\tilde T(F)$ 
over the super Teichm\"uller space $ST(F)$ of a surface $F$ with $s\geq 1$ punctures that extend the lambda length coordinates on the decorated bundle $\tilde T(F)=T(F)\times {\mathbb R}_+^s$ over the usual Teichm\"uller space $T(F)$.  In effect, the action of a Fuchsian subgroup of $PSL(2,{\mathbb R})$ on Minkowski space ${\mathbb R}^{2,1}$
is replaced by the action of a super Fuchsian subgroup of $OSp(1|2)$
on the super Minkowski space ${\mathbb R}^{2,1|2}$, where $OSp(1|2)$ denotes the orthosymplectic Lie supergroup, and the lambda lengths are extended by fermionic invariants of suitable triples of isotropic vectors in ${\mathbb R}^{2,1|2}$.  As in the bosonic case, there is the analogue of the Ptolemy transformation now on both even and odd coordinates as well as an invariant even two-form on $S\tilde T(F)$ generalizing the Weil-Petersson K\"ahler form.  This finally solves a problem posed in Yuri Ivanovitch Manin's Moscow seminar some thirty years ago to find the super analogue of decorated Teichm\"uller theory and provides a natural geometric interpretation in ${\mathbb R}^{2,1|2}$ for the super moduli of
$S\tilde T(F)$.
\end{abstract}
\maketitle

\tableofcontents

\section*{Introduction}

Let $F=F_g^s$ be a connected orientable surface of genus $g\geq 0$ with $s\geq 1$ punctures and negative Euler characteristic $2-2g-s<0$ in order that
$F={\mathcal U}/\Gamma$ is uniformized by a Fuchsian group $\Gamma$.  Namely, let ${\mathcal U}=\{ z=x+iy\in{\mathbb C}:y>0\}$ denote the upper half plane
with its Poincar\'e metric
 $ds^2=\frac{dx^2+dy^2}{y^2}$ and projective matrix group $PSL(2,{\mathbb R})=SL(2,{\mathbb R})/\pm I$ of oriented isometries, where $I$ denotes the identity matrix;
there is then an injective representation $\rho:\pi_1\to PSL(2,{\mathbb R})$ of the fundamental group $\pi_1=\pi_1(F)$, which is a free group of rank $2g+s-1$,  onto a discrete subgroup $\Gamma<PSL(2,{\mathbb R})$ so that non-trivial loops about punctures are represented by parabolic transformations, namely, those with absolute trace equal to two.  See  \cite{beardon,ford,pb} for example.

The {\it Teichm\"uller space} of $F$ is
$$T(F)~=~{\rm Hom}'(\pi_1,PSL(2,{\mathbb R}))/PSL(2,{\mathbb R}),$$
where the prime indicates Fuchsian representations as just defined and the action
of $PSL(2,{\mathbb R})$ on ${\rm Hom}'$ is by conjugation.  The {\it super-Teichm\"uller space} of $F$
as already formulated 
in the context of representation theory and moduli spaces by Bryant and Hodgkin \cite{bh,hodgkin} (see also \cite{crane-rabin}) is
$$
{ST}(F)~=~{\rm Hom}'(\pi_1, OSp(1|2))/OSp(1|2),
$$
 where the corresponding super Fuchsian representations comprising ${\rm Hom}'$ are defined to be those whose projection $\pi_1\to OSp(1|2)\to SL(2,{\mathbb R})\to PSL(2,{\mathbb R})$ are Fuchsian, where $OSp(1|2)$ denotes the { orthosymplectic group} of $(2|1)$-by-$(2|1)$ dimensional super matrices with its canonical projection $OSp(1|2)\to SL(2,{\mathbb R})$, cf. \cite{manin} or Appendix I, and the action on ${\rm Hom}'$ is again by conjugation.  The similarities are evident.  In particular, the {\it mapping class group} $MC(F)$ of homotopy classes of orientation-preserving homeomorphisms of $F$ acts on $T(F)$ and $ST(F)$ in the natural way.

\begin{figure}[hbt] 
\centering           
\includegraphics[width=0.45\textwidth]{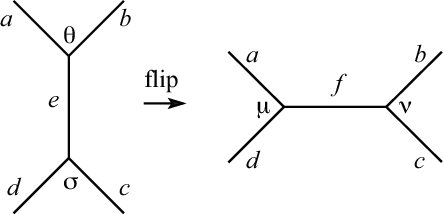}
\caption{A flip on a trivalent fatgraph with notation subsequently explained}  \label{flip0}               
\end{figure}

Consider a graph 
$\tau\subset F$ embedded in $F$ as a deformation retract also called a {\it spine} of $F$.
The {\it valence} of a vertex of $\tau$ is the number of half-edges incident upon it, where
a half-edge is defined as a complementary component to an interior point of the edge,
and $\tau$ is said to be {\it trivalent} if each vertex has valence exactly three.
An orientation on $F$ induces the counter clockwise ordering on the half edges of $\tau$ incident on each fixed vertex thus giving the abstract graph $\tau$ the structure of a 
{\it fatgraph} sometimes also called a ribbon graph.  There is a combinatorial move on
trivalent fatgraph spines $\tau\subset F$ called a {\it flip} as illustrated in Figure \ref{flip0},
where one contracts an edge of $\tau$ with distinct endpoints and then expands the resulting 4-valent vertex in the unique distinct manner in order to produce another trivalent fatgraph spine.   This leads to the so-called {\it Ptolemy groupoid} (see e.g. \cite{pb}) of $F$
whose objects are homotopy classes of trivalent fatgraph spines in $F$ and whose morphisms are compositions of flips.

As we shall recall in the next introductory section dedicated to the bosonic case, finite compositions of flips act transitively on homotopy classes of trivalent fatgraph spines.  It follows that flips generate $MC(F)$ in the sense that if $\tau\subset F$ is a trivalent fatgraph spine and $\varphi\in MC(F)$, then there is a sequence
$\varphi(\tau)=\tau_1-\tau_2-\cdots -\tau_n=\tau$ of trivalent fatgraph spines of $F$ where any
consecutive pair differ by a flip.

 In fact \cite{witten}, the components of ${ST}(F)$ are naturally indexed by the set 
 $\Omega(F)$ of spin structures on $F$, and here is a basic difference from the bosonic case: the super-Teichm\"uller space is disconnected with the $MC(F)$-action permuting components. 
 There are a number of equivalent formulations of spin structure, and we shall rely upon several of them at various junctures.
 Milnor's elegant formulation of a spin structure on $F$ is a class in the mod two first cohomology of the unit tangent bundle of $F$ which is non-zero on the fiber class;
see \cite{milnor,johnson}.  More combinatorial formulations
from the literature which we shall require are as follows:

\medskip

\leftskip .2in

\noindent $\bullet$ The
description \cite{johnson} due to Johnson in terms of {\it quadratic forms}
${\mathcal Q}(F)$ 
on $H_1=H_1(F;{\mathbb Z}_2)$,
i.e., functions $q:H_1\to{\mathbb Z}_2$
which are quadratic for the intersection pairing 
$\cdot:H_1\otimes H_1\to {\mathbb Z}_2$ in the sense that
$q(a+b)=q(a)+q(b)+a\cdot b$ if $a,b\in H_1$.

\medskip

\noindent $\bullet$ Cimasoni and Reshetikhin \cite{cr1,cr2} formulate spin structures using \cite{johnson} in terms of so-called Kastelyn orientations and dimer configurations
on the one-skeleton of a suitable CW decomposition of $F$ as we shall explain
in detail later.

\medskip

\noindent $\bullet$ A spin structure on a uniformized surface $F={\mathcal U}/\Gamma$
is determined by a lift $\tilde\rho:\pi_1\to SL(2,{\mathbb R})$ of $\rho:\pi_1\to PSL_2({\mathbb R})$, and Natanzon \cite{natanzon} computes in terms of the
quadratic form $q$ that
${\rm trace}~\tilde \rho(\gamma)>0$ if and only if $q([\gamma])\neq 0$,
where $[\gamma]\in H_1$ here and in the sequel is the image of $\gamma\in\pi_1$ under the mod two Hurewicz map.

\medskip

\leftskip=0ex

Our first main result gives yet another combinatorial formulation of
spin structures on $F$ in terms of the equivalence classes ${\mathcal O}(\tau)$ of all orientations
on a trivalent fatgraph spine $\tau\subset F$, where the equivalence relation is 
generated by reversing the orientation of each edge incident on some fixed vertex, with the added bonus of 
a computable evolution under flips.

\bigskip

\begin{figure}[hbt] 
\centering           
\includegraphics[width=0.8\textwidth]{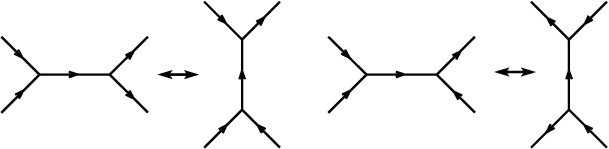}
\caption{Flip transformation for oriented fatgraphs}  \label{flipper}               
\end{figure}

\noindent{\bf Theorem A.}
{\it Fix any trivalent fatgraph spine $\tau\subset F=F^s_g$.  Then $\mathcal{O}(\tau)$ and  $\mathcal{Q}(F)$ are isomorphic as affine $H^1(F; \mathbb{Z}_2)$-spaces. Moreover, the action of $MC(F)$ on $\Omega(F)$ lifts to the action of the Ptolemy groupoid on $\mathcal{O}(\tau)$, so that under the flip transformations the  orientations of the edges change as illustrated in Figure \ref{flipper}.}\\

A spin structure on $F$ manifest for instance as a quadratic form $q\in{\mathcal Q}(F)$ distinguishes two types of punctures as follows. If $\gamma_p$ is a simple loop about the puncture $p$ of $F$, then $p$ is called a {\it Neveu-Schwarz} (NS) puncture if $q ([\gamma_p])=0$ and is otherwise called a {\it Ramond} (R) puncture.
The dimension of a component of $ST(F)$ depends on the types of punctures, and each component of ${ST(F_g^s)}$ is in fact \cite{crane-rabin},\cite{schwarz},\cite{witten},\cite{DW} a superball of dimension
$(6g-6+2s|4g-4+2n_{NS}+n_{R})$, where  $n_{NS}$ and $n_{R}$ are the respective numbers of $NS$ and $R$ punctures. It is not hard to see that $n_R$ must be even since the sum of the ${\mathbb Z}_2$-monodromies about the punctures must agree with the trivial monodromy about a small inessential curve in the surface, and we furthermore evidently have $s=n_{NS}+n_{R}$.

In the bosonic case \cite{penner,pb}, there is a principal ${\mathbb R}_+^s$-bundle $$\tilde{T}(F)=T(F)\times{\mathbb R}_+^s\to T(F)$$ called the {\it decorated Teichm\"uller space} of $F=F_g^s$ recalled in some detail in the next introductory section.   Letting ${\mathbb R}_+^{p|q}$ here and in the sequel denote the subspace of ${\mathbb R}^{p|q}$ whose even coordinates have positive body,
there is likewise a principal ${\mathbb R}_+^{s|n_R}$-bundle over each component of $ST(F)$ with $n_R$ Ramond punctures called the {\it decorated super-Teichm\"uller space} and written simply $$S\tilde T(F)\to ST(F),$$
and the action of $MC(F)$ lifts to $S\tilde T(F)$ as we shall see.
Here we combine the main results from Sections 8 and 9 into one result:

\medskip

\noindent {\bf Omnibus Theorem B.}  {\it Fix a surface $F=F_g^s$ of genus $g\geq 0$ with $s\geq 1$ punctures, where $2g-2+s>0$, let
$\tau\subset F$ be some trivalent fatgraph spine and
suppose that
$\omega$ is an orientation on the edges of $\tau$ whose class in ${\mathcal O}(\tau)$ determines the component $C$ of
$S\tilde T(F)$.  Then there are global affine coordinates on $C$, one even coordinate
called a $\lambda$-length for each edge and one odd coordinate called a $\mu$-invariant for each vertex of $\tau$,
 the latter of which are taken modulo an overall change of sign denoted by ${\mathbb Z}_2$; that is, $\lambda$-lengths and $\mu$-invariants
establish a real-analytic homeomorphism 
$$C\to{\mathbb R}_+^{6g-6+3s|4g-4+2s}/{\mathbb Z}_2.$$
These coordinates are natural in the sense that
if $\varphi\in MC(F)$ has induced action $\tilde\varphi$ on $\tilde\Gamma\in S\tilde T(F)$, then
$\tilde\varphi(\tilde\Gamma)$ is determined 
by the orientation and coordinates on edges and vertices of $\varphi(\tau)$ 
induced by $\varphi$ from the orientation $\omega$, the $\lambda$-lengths and $\mu$-invariants on $\tau$.\\

\noindent Orientations on fatgraph spines evolve under flips in accordance with the previous theorem, and the following
{\rm super Ptolemy transformation} further describes the evolution of $\lambda$-lengths and $\mu$-invariants in the notation of Figure \ref{flip0} where nearby
Roman letters denote $\lambda$-lengths, nearby Greek letters denote $\mu$-invariants
and  $Z=\frac{ac}{bd}$ denotes the {\rm cross-ratio}
\begin{eqnarray}
&& ef=(ac+bd)\Big(1+\frac{\sigma\theta\sqrt{Z}}{1+Z}\Big),\nonumber\\
&&\nu=\frac{\sigma+\theta\sqrt{Z}}{\sqrt{1+Z}},\quad
\mu=\frac{\sigma\sqrt{Z}-\theta}{\sqrt{1+Z}}.\nonumber
\end{eqnarray} \\

\noindent Finally, there is an even 2-form on $S\tilde T(F)$
which is invariant under super Ptolemy transformations, namely,
\begin{eqnarray}\label{wpsuper}
\omega=
\sum_{v}d\log a\wedge d\log b+d\log b\wedge d\log c
+d\log c\wedge d\log a -(d\theta)^2,\nonumber
\end{eqnarray}
where the sum is over all vertices $v$ of $\tau$ where the consecutive half edges incident on $v$ in clockwise order have induced $\lambda$-lengths $a, b, c$ and $\theta$ is the $\mu$-invariant of $v$.
}
\medskip

\noindent In order to further explain the coordinates and give an intrinsic meaning to
the decorated spaces, we shall recall and extend the bosonic case in the next 
section.

The super-Teichm\"uller space is of interest on its own as one of the higher Teichm\"uller theories associated with the simplest supergroup extension of $PSL(2,\mathbb{R})$, namely, the orthosymplectic group $OSp(1|2)$ which is however special among supergroups in that its invariant bilinear form is non-degenerate. Nevertheless, this should be the first step in a general approach to higher super Teichm\"uller theory
as well as its quantization.  Notice that the symplectic or corresponding Poisson structure of Theorem B has constant coefficients opening the possibility of canonical quantization  as in the bosonic case \cite{Fock-Chekhov,Kashaev}. 

Furthermore just as for $\lambda$-lengths in the bosonic case, our coordinates on decorated super-Teichm\"uller space provide a computationally effective description of super moduli geometry. 
Another key reason for interest in the super-Teichm\"uller space is that it is a cornerstone of superstring perturbation theory, and the geometry of the supersymmetric moduli uncovered here is evidently more involved than its bosonic counterpart.  
Moreover, the well-known relationship \cite{verlinde} of Teichm\"uller theory with $(2+1)$-dimensional gravity may suggest that the super-Teichm\"uller theory as probed here may play an analogous role for $(2+1)$-dimensional supergravity.  

We finally mention
that prior to our work, there was a PhD thesis \cite{zababaher}, to be continued in \cite{BCF}, where the author provided 
other coordinates on a version of super-Teichm\"uller space. The author used a different combinatorial description of spin structures, by marking sectors in the triangulation, and his contruction is based upon quite a different approach (known as path-ordered method, see e.g. \cite{pb}) effectively using so-called shear coordinates instead of $\lambda$-lengths without our connection to super Minkowski geometry described in the next several sections.\\

\noindent{\bf Acknowledgements.} We are grateful to Institut des Hautes \'Etudes Scientifiques for hospitality and wonderful atmosphere at the initial phase of this work. We would like to thank R. Donagi, R. Kaufmann, Yu.I. Manin, J. Rabin, A.S. Schwarz, A. Voronov, E. Witten for their interest, comments and discussions.  A.M.Z. 
is indebted to M. Liu, M. Khovanov and A. Okounkov for their encouragement and support during his stay at Columbia University and to A.N. Fedorova for careful reading of the manuscript.

\section*{Bosonic background and beyond} 

The decorated Teichm\"uller space $\tilde T(F)$ is intimately connected to the geometry of Minkowski space ${\mathbb R}^{2,1}$, namely, ${\mathbb R}^3$ imbued with the 
quadratic form $z^2-x^2-y^2$ and its corresponding bilinear pairing
$ <(x,y,z),(x',y',z')>=zz'-xx'-yy'$ or equivalently the quadratic form
$x_1x_2-y^2$ in the variables $x_1=z-x$ and $x_2=z+x$.  The upper sheet  
${\mathbb H}=\{ (x,y,z)\in{\mathbb R}^3:z^2-x^2-y^2=1~{\rm and}~z> 0\}$ of the hyperboloid
with the induced metric gives a model for the
hyperbolic plane \cite{beardon}.  Furthermore \cite{penner,pb},  the positive light cone
${L}^+=\{ (x,y,z)\in{\mathbb R}^3:z^2-x^2-y^2=0~{\rm and}~z> 0\}$
parametrizes horocycles (that is, those curves with 
geodesic curvature unity)  in the sense that $h(u)=\{ w\in{\mathbb H}:<u,w>=1\}$ establishes a bijection\footnote{
In fact, any positive constant will suffice here, and $2^{-{1\over 2}}$ is the more natural choice; cf.\ \cite{penner,pb}.} between all $u\in L^+$ and the collection of all horocycles $h(u)\subset{\mathbb H}$, and moreover, this identification is geometrically
natural in the sense that ${1\over 2} \log{<u,v>}$ is  the signed hyperbolic distance between $h(u)$ and $h(v)$.  This invariant $\sqrt{<u,v>}$ of a pair $\{h(u),h(v)\}$ of horocycles  is called the {\it lambda length}, and 
these are the basic coordinates on $\tilde T(F)$.

It is convenient both here and in the sequel to consider not only trivalent fatgraph spines of $F=F_g^s$ but also their duals:  an {\it ideal triangulation} $\Delta$ of $F$ is a maximal family of arcs embedded in $F$ with endpoints at the punctures, which are here regarded as distinguished points of $F$, where no two arcs in $\Delta$  are properly homotopic
or intersect except perhaps at their endpoints.
By maximality, each complementary region to $\Delta$ in $F$ is an ideal triangle, and these meet along their frontiers in $F$.  Construct a trivalent fatgraph spine $\tau=\tau(\Delta)$ of $F$ with one vertex for each complementary region and one edge for each arc in $\Delta$ required to connect the vertices corresponding to regions on either side of the arc.  This construction evidently establishes a bijection between homotopy classes of trivalent fatgraph spines and homotopy classes of ideal triangulations of $F$ together with a natural identification between their edges.  The dual of a flip, also called a flip, is the removal of an edge from $\Delta$ with distinct triangles on either side followed by its replacement
by the other diagonal of the quadrilateral formed by their union.

A key point of the decorated space in the bosonic case is that the fiber  ${\mathbb R}^s_+$ of $\tilde T(F)=T(F)\times {\mathbb R}^s_+$ is identified with all $s$-tuples of (lengths of) horocycles in $F$ with one horocycle about each puncture; this is precisely the sense in which a usual hyperbolic structure in $T(F)$ is {\sl decorated
with a horocycle at each puncture} in $\tilde T(F)$.  
The $MC(F)$-action on $T(F)$ thus lifts to $\tilde T(F)$ by permuting lengths.
An arc connecting punctures in a decorated hyperbolic surface has a well-defined lambda length computable in the surface itself or equivalently in ${\mathbb R}^{2,1}$
as just discussed.  \\

\noindent {\bf Omnibus Theorem C.}  \cite{penner,penner-jdg2,pb} {\it Fix a surface $F=F_g^s$ of genus $g\geq 0$ with $s\geq 1$ punctures, where $2g-2+s>0$, and let
$\Delta$ be a homotopy class of ideal triangulation  or equivalently of a trivalent fatgraph spine of $F$.  Then the assignment of lambda lengths to the arcs in $\Delta$ establishes a real-analytic homeomorphism $\tilde T(F)\to{\mathbb R}_+^\Delta$.  
Moreover, the lambda length $\lambda(\alpha;\tilde\Gamma)$ of a (homotopy class of) arc $\alpha$ in $F$ connecting punctures for $\tilde \Gamma\in \tilde T(F)$ is natural in the sense that
if $\varphi\in MC(F)$ has induced action $\tilde \varphi$ on $\tilde T(F)$, then 
$\lambda(\alpha;\tilde\Gamma)=\lambda(\varphi(\alpha);\tilde\varphi(\tilde\Gamma))$.\\

\noindent The Ptolemy transformation $ef=ac+bd$ describes the evolution of lambda lengths under flips in the notation of Figure \ref{flip0}.\\

\noindent The Weil-Petersson K\"ahler form on $M(F)$ pulls back to the 
Ptolemy-invariant form
\begin{eqnarray}
\omega=
2\sum_{v}d\log a\wedge d\log b+d\log b\wedge d\log c
+d\log c\wedge d\log a,\nonumber
\end{eqnarray}
where the sum is over all complementary triangles to $\Delta$ with consecutive half edges in clockwise order having induced lambda lengths $a, b, c$.\\

\noindent A convex hull construction in ${\mathbb R}^{2|1}$ gives rise to a real-analytic  $MC(F)$-invariant ideal cell decomposition of $\tilde T(F)/{\mathbb R}_+$ itself where there is one open simplex together with certain of its faces for each homotopy class of
decompositions of $F$ into ideal polygons and the face relation is generated by removal of arcs.}\\

Thus in the equivalent formalism of trivalent fatgraph spines as opposed to ideal triangulations, Theorem B extends all but the last paragraph of Theorem C 
from $\tilde T(F)$ to $S\tilde T(F)$, and the proof of the former provides a paradigm for the proof of the latter.  Notice that by this last paragraph of Theorem C, we may connect interior points of simplices for any two ideal triangulations of $F$ by a path in $\tilde T(F)/{\mathbb R}_+$, which 
can be perturbed to general position with respect to the codimension-one faces
of the cell decomposition.  Since crossing these faces corresponds
to flips, it follows that any two (homotopy classes of) ideal triangulations of $F$ are related by a finite sequence of flips, and dually we recover Whitehead's result \cite{flp}:

\bigskip

\noindent{\bf Corollary D} {\it For any surface $F=F_g^s$ with $2g-2+s>0$, finite sequences of flips act transitively on homotopy classes of trivalent fatgraph spines in $F$.
}

\bigskip

\noindent Moreover,
the codimension-two faces analogously give rise to a presentation of the Ptolemy groupoid of $F$, cf. \cite{penner,pb}.

We next discuss the proof of the first part of Theorem C, namely, the construction of lambda length coordinates on $\tilde T(F)$ which brings us more deeply into Minkowski space.   The topological universal cover $\tilde F$ of $F$ may be identified with upper half space ${\mathcal U}$ or equivalently via the Cayley transform
$z\mapsto \frac{z-i}{z+1}$ with the unit disk ${\mathbb D}$ supporting its Poincar\'e metric and ideal boundary the circle ${\mathbb S}^1$ at infinity.  Central projection of ${\mathbb H}$ from
$(0,0,-1)\in {\mathbb R}^{2,1}$ to the disk at height zero establishes an isometry of
${\mathbb H}$ and ${\mathbb D}$ which continuously extends to the projection $L^+\to {\mathbb S}^1$ mapping $u\in L^+$ to the center of the horocycle $h(u)$.

An ideal triangulation $\Delta$ of $F$ lifts to an ideal triangulation $\tilde\Delta$ of $\tilde F$, and the collection of ideal points of $\tilde\Delta_\infty\subset {\mathbb S}^1$ is invariant under homotopy of $\Delta$ in $F$.  In order to define a point of $\tilde T(F)$, we must determine
a Fuchsian representation $\rho:\pi_1\to SO_+(2,1)\approx PSL(2,{\mathbb R})$ in the component $SO_+(2,1)$ of the identity of the Minkowski isometry group $SO(2,1)$, 
corresponding to the underlying point in $T(F)$, together with a lift
$\ell:\tilde \Delta_\infty\to L^+$, corresponding via affine duality to the decoration and realizing the lambda lengths in the obvious sense, which is $\pi_1$-equivariant with respect to our constructed representation $\rho:\pi_1\to SO_+(1,2)$, namely, we have
$\ell(\gamma(p))=\rho(\gamma) (\ell(p))$, for all $p\in\tilde\Delta_\infty$ and $\gamma\in\pi_1$.   This construction of $\ell$ and $\rho$ from lambda lengths is performed recursively as we shall recall in proving Theorem B.

In our current case, there is an embedding of $OSp(1|2)$ into the super Lorentz group of the super Minkowski space ${\mathbb R}^{2,1|2}$ with pairing
$x_1x_2-y^2+2\phi\theta$ as described in the next section.
There is again a positive light cone $\hat L^+$ comprised of isotropic vectors so that the bodies of $x_1$ and $x_2$ are non negative.  However, here is another fundamental distinction between the bosonic case and the general case treated here: whereas the action of $SO_+(2,1)$ on $L^+$ is transitive, an $OSp(1|2)$-orbit
of positive isotropic vectors is determined by a fermionic invariant $\pm\xi$ up to sign.  Though most of the computations of this paper can be completed in the general setting, the super Teichm\"uller theory seems to require taking this fermion label $\pm \xi=0$.  It is the special light cone $\hat L^+_0\subset \hat L^+$  consisting of those positive isotropic vectors with vanishing fermion label $\xi=0$ that provide the analogue of $L^+$ for us here.  Again, we define the $\lambda$-length of a pair of points in $\hat L^+_0$ to be the square root of their inner product and prove that this is the unique invariant of the $OSp(1|2)$-orbit of a pair of linearly independent points in $\hat L^+_0$.

We shall again recursively define a mapping $\ell:\tilde\Delta_\infty\to \hat L^+_0$ which realizes $\lambda$-lengths in the obvious sense and is $\pi_1$-equivariant with respect to the representation $\hat\rho:\pi_1\to OSp(1|2)$ that we construct.  Here is yet another fundamental distinction between the bosonic case and the general case: whereas  $SO_+(2,1)$ acts transitively on triples of rays in $L^+$ which are consistent with the positive orientation on ${\mathbb R}^{2,1}$, an $OSp(1|2)$ orbit of a triple in 
$\hat L^+_0$ whose underlying bosonic vectors in $L^+$ have this property is again determined by a fermion invariant $\pm\mu$ up to sign.  Manin \cite{manin} had already observed this basic phenomenon hence our term $\mu$-invariant for the odd invariants associated to vertices of a fatgraph spine in Theorem B which come from consistent choices of signs and are, like $\lambda$-lengths, realized by the mapping $\ell$ in the obvious sense.  Consistency here is given by an explicit relationship on signs of $\mu$-invariants for adjacent triangles in $\tilde \Delta$.

There are thus three basic differences here from the bosonic case: the failure of transitivity of the $OSp(1|2)$-action on points and on triples in $\hat L^+_0$ already mentioned and the further fact that the identification of $L^+$ with the space of horocycles in ${\mathbb H}$ has no known analogue in the general case.  Thus, the decorated super-Teichm\"uller space $S\tilde T(F)$ {\sl can only be defined here} as the space of $OSp(1|2)$ orbits of those maps 
$\ell:\tilde\Delta_\infty\to \hat L^+_0$ that are $\pi_1$-equivariant for some super Fuchsian representation with no intrinsic interpretation of {\sl super horocycle}  for the decoration beyond the analogous but un-illuminating affine dual in ${\mathbb R}^{2,1|2}$ of a point in $\hat L^+_0$.  Indeed, in addition to the research frontiers discussed at the end of the previous section that the current work presumably illuminates, so also first glimpses of {\sl super hyperbolic geometry} are hopefully to be gleaned here.

\section{The hyperboloid, light cone and $OSp(1|2)$-action}
The supergroup $OSp(1|2)$ is defined as follows. The group elements are square $(2|1)\times(2|1)$ supermatrices with superdeterminant equal to 1 which satisfy the relation
\begin{eqnarray}\label{ospdef}
g^{st}Jg=J,
\end{eqnarray}
where the superscript $st$ denotes the super transpose and $J= \bigl(\begin{smallmatrix}
\hskip .6em0&1&0\\ -1&0&0\\ \hskip .6em0&0&1\\
\end{smallmatrix} \bigr)$. We refer the reader to Appendix~I for more information about $OSp(1|2)$ 
including the definition of super transpose (or see \rf{action} below) and
our sign conventions for products of supermatrices.
The useful property 
\begin{eqnarray}\label{ainv}
g^{st}=Jg^{-1}J^{-1}
\end{eqnarray}
is of course a direct consequence of  \rf{ospdef}.
We are interested in the adjoint action of $OSp(1|2)$ and to this end consider its even element 
\begin{eqnarray}
N_0=-yh+x_1 X_{-}-x_{2}X_{+}+\phi v_{-}-\theta v_{+}.
\end{eqnarray}
We claim that the adjoint action $N_0\mapsto g^{-1}N_0g$ of $g\in OSp(1|2)$ on $N_0$ leaves invariant $x_1x_2-y^2+2\phi\theta$
since it is proportional to the quadratic form arising from the Killing form of $OSp(1|2)$ applied to $N_0$. 

We shall prove this differently as follows. Observe that in the 3-dimensional representation of $OSp(1|2)$ owing to the property \rf{ainv}, 
the element $M_c=JN_c$, where $N_c=N_0+cI$ for any fixed constant $c$, transforms as 
\begin{eqnarray}
N_c\mapsto g^{-1}N_cg \quad\Rightarrow \quad M_c\mapsto g^{st}M_cg
\end{eqnarray}
under the adjoint action. 
In particular when $c$ is invertible, this implies that the 
superdeterminant of 
\begin{eqnarray}\label{goop}
M_c=\left( \begin{array}{ccc}
x_1 & y-c & \phi \\
y+c & x_2 & \theta \\
-\phi & -\theta & c \end{array} \right)
\end{eqnarray}   
is preserved under the action of $OSp(1|2)$ sending $M_c$ to $g^{st}M_cg$. It is not hard to calculate that
\begin{eqnarray}
sdet(M_c)=\frac{1}{c}(x_1x_2-y^2+2\phi\theta+c^2),
\end{eqnarray}  
and so 
$x_1x_2-y^2+2\phi\theta$ is invariant under the action of $OSp(1|2)$ as was claimed. 
The following proposition therefore holds.
\begin{Prop}
The formula $M_c\mapsto g^{st}M_cg$, for fixed but arbitrary $c$, gives the action of the $OSp(1|2)$ subgroup of the full Lorentz  supergroup
of the superspace $\mathbb{R}^{2,1|2}$ with the Minkowski pairing defined by the quadratic form $x_1x_2-y^2+2\phi\theta$.
\end{Prop}   

\noindent To be entirely explicit, the pairing of two vectors $A=(x_1,x_2, y,\phi, \theta)$ and $A'=(x'_1,x'_2, y',\phi', \theta')$ in $\mathbb{R}^{2,1|2}$ is given by
\begin{eqnarray}
\langle A, A'\rangle=\frac{1}{2}(x_1x'_2+x'_1x_2)-yy'+\phi\theta'+\phi'\theta.
\end{eqnarray} 
In keeping with \cite{penner,pb}, we shall henceforth refer to the square root of such an inner product as a {\it $\lambda$-length}.

Two surfaces of special importance for us in the following are the 
{\it (super) hyperboloid} $\hat{\mathbb{H}}$ consisting of points $A\in\mathbb{R}^{2,1|2}$ satisfying 
the condition $\langle A, A\rangle=1$ corresponding to $c=1$ in equation \rf{goop}, where the bodies of the $x_1$- and $x_2$-coordinates of $A$ are non-negative, 
and most especially the {\it (positive super) light cone} $\hat{L}^+$
consisting of points $B\in \mathbb{R}^{2,1|2}$ satisfying $\langle B, B\rangle=0$ and corresponding to $c=0$, where again the bodies of $x_1$- and $x_2$-coordinates 
are non-negative.

A standard superspace which however plays a subsidiary role for us, the complex {\it superplane} $\mathbb{C}^{1|1}$ consists of pairs $(z,\eta)$, and its subspace the {\it super upper half-plane} $\hat{\mathcal{U}}$ is comprised of those points $(z,\eta)$ such 
that the body of the real part of $z$ is non-negative. 
It is well known \cite{crane-rabin,schwarz, witten} that $OSp(1|2)$ acts on $\hat{\mathcal{U}}$ by means of superconformal transformations
\begin{eqnarray}\label{transff}
&&z\to \frac{az+b}{cz+d}+\eta\frac{\gamma z+\delta}{(cz+d)^2}, \nonumber\\
&&\eta\to \frac{\gamma z+\delta}{cz+d}+\eta\frac{1+\2 \delta\gamma}{cz+d}, 
\end{eqnarray}

Another direct analogue of the standard bosonic case, we have
\begin{Thm}\label{splane}
The expressions
\begin{eqnarray}
\eta=\frac{\theta}{x_2}(1+iy)-i\phi,\quad z=\frac{i-y-i\phi\theta}{x_2} 
\end{eqnarray}
define an $OSp(1|2)$-equivariant monomorphism from the hyperboloid $\hat{\mathbb H}$ onto the super half-plane $\hat{\mathcal{U}}$.
\end{Thm}

\noindent{\bf Proof.} The easiest way to prove this statement is to consider the infinitesimal actions of the corresponding generators
described in Appendix I. For example, the transformation $M_c\mapsto (\exp(\alpha v_+))^{st} M_c \exp(\alpha v_+)$ amounts to the infinitesimal
\begin{eqnarray}
&&\hskip 1ex\delta_{\alpha}y=-\alpha\phi, \quad\hskip.5ex \delta_{\alpha}\phi=\alpha x_1, \quad \delta_{\alpha}x_1=0,\nonumber\\
&&\delta_{\alpha}x_2=-2\alpha\theta, \quad \delta_{\alpha} \theta=\alpha y,
\end{eqnarray} 
and therefore
\begin{eqnarray}
\delta_{\alpha}z=\alpha\biggl (\frac{2\theta}{x_2^2}(i-y)+\frac{1}{x_2}(\phi-ix_1\theta+i\phi y)\biggr ).
\end{eqnarray} 
Meanwhile, we have
\begin{eqnarray}
&&\eta z=i\phi\frac{i-y}{x_2}-i\theta\frac{(y-i)^2}{x_2^2}\nonumber\\
&&\hskip .2in=\frac{i\phi}{x_2}+\frac{\phi}{x_2}-\frac{i\theta y^2}{x_2^2}-\frac{2\theta y}{x_2^2}+\frac{i\theta}{x_2^2}\nonumber\\
&&\hskip .2in=\frac{i\phi y}{x_2}+\frac{\phi}{x_2}-\frac{ix_1}{x_2}-\frac{2\theta y}{x_2^2}+\frac{2i\theta}{x_2^2},
\end{eqnarray}
where we have used the relation $x_1x_2-y^2+2\phi\theta=1$ in the last line.  It follows that $\delta_{\alpha}z=\alpha\eta z$
as required, and one can similarly show that $\delta_{\alpha}\eta=-\alpha z$. Thus, this corresponds to the superconformal transformation. We leave it for the reader to complete the proof for the other four infinitesimal transformations corresponding to the generators discussed in Appendix I.
\hfill $\blacksquare$

\medskip

Again just to be entirely explicit in the context of relevant subsequent calculations, the light cone is described by
$$\hat{L}^+=\biggl \{ \bigl(\begin{smallmatrix}
~x_1&~y&\phi\\ ~y&~x_2&\theta\\-\phi&-\theta&0\\
\end{smallmatrix} \bigr):x_1x_2-y^2+2\phi\theta=0~{\rm and}~x_1,x_2~{\rm have~non-negative~bodies}
\biggr \},$$
and the action of the supermatrix $g=\bigl(\begin{smallmatrix}
a&b&\alpha\\ c&d&\beta\\ \gamma&\delta&f\\
\end{smallmatrix} \bigr)\in OSp(1|2)$ on $A\in\hat{L}^+$ is given by
\begin{eqnarray}\label{action}
g\cdot A=g^{st} A g= \biggl(\begin{matrix}
~a&~c&\gamma\\ ~b&~d&\delta\\ -\alpha&-\beta&f\\
\end{matrix} ~\biggr) ~~
\biggl(\begin{matrix}
~x_1&~y&\phi\\ ~y&~x_2&\theta\\-\phi&-\theta&0\\
\end{matrix}~ \biggr) ~~\biggl(~\begin{matrix}
a&b&\alpha\\ c&d&\beta\\ \gamma&\delta&f\\
\end{matrix} ~\biggr),
\end{eqnarray}
where these products of supermatrices have the signs\footnote{
Namely,
$\biggl ( \begin{smallmatrix}
a_1&b_1&\alpha_1\\c_1&d_1&\beta_1\\\gamma_1&\delta_1&f_1\\
\end{smallmatrix}\biggr )
\biggl ( \begin{smallmatrix}
a_2&b_2&\alpha_2\\c_2&d_2&\beta_2\\\gamma_2&\delta_2&f_2\\
\end{smallmatrix}\biggr )=
\biggl (\begin{smallmatrix}
a_1a_2+b_1c_2-\alpha_1\gamma_2&a_1b_2+b_1d_2-\alpha_1\delta_2&a_1\alpha_2+b_1\beta_2+\alpha_1 f_2\\
c_1a_2+d_1c_2-\beta_1\gamma_2&c_1b_2+d_1d_2-\beta_1\delta_2&c_1\alpha_2+d_1\beta_2+\beta_1 f_2\\
\gamma_1 a_2+\delta_1 c_2+\delta_1\gamma_2&\gamma_1 b_2+\delta_1 d_2+f_1 \delta_2&-\gamma_1\alpha_2-\delta_1\beta_2+f_1f_2\\
\end{smallmatrix}\biggr ).
$}
explained in Appendix I.
This is entirely analogous to the bosonic action  \cite{penner,pb}  of
$PSL(2,{\mathbb R})\approx SO_+(2,1)$ on  the light cone in ${\mathbb R}^{2,1}$ given by the change of basis for binary symmetric bilinear forms.

There is the particular element of $OSp(1|2)$ given by $g^r=\bigl(\begin{smallmatrix}
-1&\hskip .6em 0&~0\\ \hskip .6em 0&-1&~0 \\ \hskip .6em 0&\hskip .6em 0&~1\\
\end{smallmatrix} \bigr)$ that is of special significance.  The supermatrix $g^r$ generates the center of $OSp(1|2)$, and its explicit action on any
$A=(x_1,x_2,y,\phi,\theta)$ is given by  $g^r\cdot A=(x_1,x_2,y,-\phi,-\theta)$. Thus, $g^r$ simply changes the signs of the
fermions and will henceforth be referred to as {\it (fermionic)~reflection}.

\section{Orbits of $OSp(1|2)$ in the light cone}
We next show that $OSp(1|2)$ does not act transitively on the light cone $\hat{L}^+$, and
in fact, the moduli space of orbits is homeomorphic to the space $\mathbb{R}^{0|1}/{\mathbb Z}_2$, 
where ${\mathbb Z}_2$ acts by the change of sign of fermions.  To begin, we normalize
with respect to the subgroup $SL(2,{\mathbb R})<OSp(1|2)$.

\begin{Lem}
For each vector $A\in \hat{L}^+$, there is some  $g\in SL(2,{\mathbb R})<OSp(1|2)$ so that
$g\cdot A=t(1, 1, 1+\phi\psi, \phi,\psi)$, where t has positive body.
\end{Lem}
\noindent {\bf Proof.} Consider an arbitrary vector $(x_1,x_2, y, \rho, \lambda)\in \hat{L}^+$. Since one of $x_1$ or $x_2$ is invertible,
we can apply an element of the $SL(2,{\mathbb R})$ subgroup to transform to a vector $(x'_1,x'_2, y ', \rho ', \lambda ')$ where both of 
$x_1', x_2'$ are invertible and the body of $y'$ is positive. We can subsequently apply a diagonal 
matrix from the $SL(2,{\mathbb R})$ subgroup in order that the resulting vector $(x''_1,x''_2, y'', \rho '', \lambda '')$ satisfies $x''_1=x''_2$
and hence has the required form.
\hfill $\blacksquare$\\

The next result provides the classification of orbits, namely,  we can reduce all degrees of freedom to a single fermion modulo sign
via the action of $OSp(1|2)$.

\begin{Prop}\label{orb}
Every vector in the light cone can be put into the form $$e_{\theta}=(1,0,0,0,\theta)\in\hat{L}^+$$ via an $OSp(1|2)$ transformation. The only solutions to the equation $e_{\theta'}=g\cdot e_{\theta}$, where $g\in OSp(1|2)$, are given by $\theta'=\pm \theta$.
\end{Prop}

\noindent {\bf Proof.} According to the previous lemma, in order to prove the first part, we may assume that our
specified vector is of the form 
$A^t_{\phi\psi}=t(1, 1, 1+\phi\psi, \phi,\psi)$, where $t$ has positive body.
A direct computation then shows that the matrix 
\begin{eqnarray}\label{gtpp}
 g^t_{\phi,\psi}=\left( \begin{array}{ccc}
0 & -\sqrt{t} & 0 \\
\frac{1}{\sqrt{t}} & \sqrt{t}(1+\phi\psi) & -\psi \\
0 & \sqrt{t}\psi & 1 \end{array} \right)
\end{eqnarray}
achieves the required expression $g^t_{\phi,\psi}\cdot A^t_{\phi\psi}=e_{\theta}$, where in fact
$\theta=t\sqrt{t}(\psi-\phi)$ and $\sqrt{t}$ is likewise taken with positive body.

The second part is proven by explicitly solving the equation
$e_{\theta'}=g\cdot e_{\theta}$ as follows.
Consider an arbitrary element
\begin{eqnarray}
g=\left( \begin{array}{ccc}
a & b & \alpha \\
c & d & \beta \\
\gamma & \delta & f \end{array} \right) \in OSp(1|2).
\end{eqnarray}
The vector $g\cdot e_\theta=(x_1,x_2, y, \rho, \lambda)$ is characterized by the identities
\begin{eqnarray}\label{ga}
&&x_1=a^2+2\gamma\theta c, \quad \rho=a\alpha+\gamma \theta\beta+c\theta f,\nonumber\\
&&x_2=b^2+2\delta\theta d, \quad \lambda=b\alpha+\delta\theta\beta+d\theta f,\nonumber\\
&&y=ab+\gamma\theta d-c\theta\delta.
\end{eqnarray}
Thus, if $g\cdot e_\theta=e_{\theta'}$, then using $y=0$ as well as the constraints
on the entries of $g$ given in Appendix I, we find that
$\rho =0$, $\lambda=\theta'$ imply
\begin{eqnarray}
\alpha=-\frac{c}{a}\theta,\quad a\theta'=\theta,
\end{eqnarray}
and then $x_1=1$, $x_2=0$ imply
\begin{eqnarray}
a=\pm(1\mp c\beta\theta), \quad b=\pm \theta \beta.
\end{eqnarray}
It follows that $\theta'=\pm \theta$ as was claimed.\hfill $\blacksquare$\\

The next result follows immediately.

\begin{Cor}  
The moduli space of orbits of the $OSp(1|2)$ action on the light cone is given by $\mathbb{R}^{0|1}/{\mathbb Z}_2,$
where ${\mathbb Z}_2$ reflects the sign of the fermion.
\end{Cor}

The explicit solution to the equation $e_{\theta}=g\cdot e_{\theta}$  given in Proposition \ref{orb} yields the following Corollary 
which will be of utility in the sequel.
\begin{Cor}\label{stab}
For $\theta\neq 0$, an element $g^s$ of the stabilizer subgroup of $e_{\theta}$ in $OSp(1|2)$ necessarily has the form
\begin{eqnarray}\label{stabil}
g^s=\left( \begin{array}{ccc}
1+c\theta\beta& \theta \beta &  -c\theta \\
 c & 1+c\theta\beta & \beta \\
\beta+c^2\theta & c\theta & 1+c\beta\theta \end{array} \right),
\end{eqnarray}
where $c, \beta$ are free parameters with $c$ even and $\beta$ odd.
Moreover for $\theta=0$, the stabilizer of $e_0$ has one component for Ramond punctures given by  the expression $g^s=\biggl(\begin{smallmatrix}
1&0&0\\c&1&\beta\\\beta&0&1\\\end{smallmatrix}\biggr)$ above as well as a second component for Neveu-Schwarz punctures given by its composition
$g^rg^s=\biggl(\begin{smallmatrix}
-1&\hskip.5em0&\hskip .6em 0\\-c&-1&-\beta\\\hskip .5em \beta&\hskip .5em 0&\hskip .5em 1\\\end{smallmatrix}\biggr)$ with the fermionic reflection $g^r$.
\end{Cor}


Given a point $A\in\hat{L}^+$, the fermion $\pm\theta$ (defined up to an overall sign) so that $A$ and $e_\theta$ lie in the same $OSp(1|2)$-orbit
is called the {\it fermion~label} of $A$, and it admits the following simple expression.

\begin{Prop}\label{orbitform}
If $A=(x_1,x_2, y, \rho,\lambda)\in \hat{L}^+$, then the fermion label 
of $A$ is given by 
$\pm\theta=\sqrt{x_1}\lambda-\frac{y}{\sqrt{x_1}}\rho$ if $x_1$ is invertible and by
$\pm\theta=\sqrt{x_2}\rho-\frac{y}{\sqrt{x_2}}\lambda$ if $x_2$ is invertible. 
\end{Prop}

\noindent Here and throughout since $x_1$- and $x_2$-coordinates have non-negative body on the positive light cone, 
there are well-defined square roots $\sqrt{x_1}$ and $\sqrt{x_2}$ also with non-negative body.  
The fermionic reflection on $\lambda,\rho$ thus changes the sign of $\theta$ here.\\

\noindent {\bf Proof.} The result follows from direct calculation starting from the formulas \rf{ga}
using the constraints on entries of $OSp(1|2)$ in Appendix I
as we explicate  in the case where $x_1$ is invertible.  We have
\begin{eqnarray}
&&\sqrt{x_1}\lambda=(a+\frac{\gamma\theta c}{a})(b\alpha+\delta\theta\beta+d\theta f) \nonumber\\
&&\hskip .38in = ab\alpha-ad \alpha\theta\beta+ad\theta f+cb\theta\alpha\beta.
\end{eqnarray}
At the same time, we have
\begin{eqnarray}
&&\frac{y}{\sqrt x_1}=
(ab+\gamma\theta d-c\theta\delta)
\frac{1}{a}(1-\frac{\gamma\theta c}{a^2})=b+\frac{\gamma\theta}{a^2}-\frac{c\theta\delta}{a},
\end{eqnarray}
and therefore
\begin{eqnarray}
&&\frac{y}{\sqrt{x_1}}\rho=(b+\frac{\gamma\theta}{a^2}-\frac{c\theta\delta}{a})(a\alpha+\gamma\theta\beta+c\theta f)\nonumber\\
&&\hskip .45in = ab\alpha-cb\alpha\theta\beta+bc\theta f - ad \theta \beta\alpha.
\end{eqnarray}
It follows that 
\begin{eqnarray}
\sqrt{x_1}\lambda-\frac{y}{\sqrt{x_1}}\rho=(ad-bc)\theta f=\theta,
\end{eqnarray}
where we again use the constraints on elements of  $OSp(1|2)$  from Appendix I.
\hfill $\blacksquare$\\

The $OSp(1|2)$-orbit of $e_0\in\hat{L}^+$ will play a special role for us here.  We shall denote it
$\hat{L}_0^+=OSp(1,2)\cdot e_0\subset\hat{L}^+$ and refer to it as the {\it special light cone}.

\begin{Cor}
The special light cone $\hat{L}_0^+$ is isomorphic to superprojective space $\mathbb{R}P^{1|1}$ with the action of $OSp(1|2)$ given by superconformal transformations. Provided $x_2\neq 0$, the natural correspondence between $(x_1,x_2,y,\phi,\psi)\in\hat{L}_0^+$ and $(z, \eta)\in \mathbb{R}P^{1|1}$ is given by:
\begin{eqnarray}
z=\frac{-y}{x_2}, \quad \eta=\frac{\psi}{x_2}.
\end{eqnarray}
\end{Cor}
\noindent The proof follows along the lines of Theorem \ref{splane}.

\section{Orbits of isotropic independent pairs and positive triples}

This section provides abstract identifications for the spaces of $OSp(1|2)$-orbits of linearly independent ordered pairs and certain triples of points
in the special light cone $\hat{L}^+_0$.  Specifically, we shall say that an ordered triple $ABC$ of points in $\hat{L}^+_0$ is {\it positive} provided $A,B,C$ are linearly independent and the underlying triple of bosonic vectors in the usual light cone $L^+\subset{\mathbb R}^{2,1}$ in this order provides a positively oriented basis
for ${\mathbb R}^{3}$ with its usual orientation.  In fact, the latter moduli space of positive triples in the light cone plays a key role in the sequel, and several further parametrizations of it are derived in the next section.

The moduli space of all ordered pairs  of linearly independent vectors  (i.e., with non-vanishing Minkowski inner product or equivalently vectors  lying in distinct rays)
in the special light cone is described by

\begin{Lem}\label{2pts}
There is a unique $OSp(1|2)$-invariant of two linearly independent vectors $A,B\in \hat{L}^+_0$, and it is given by the pairing $\langle A,B\rangle $.
\end{Lem}
\noindent {\bf Proof.} According to Proposition \ref{orb}, by applying an appropriate element of $OSp(1|2)$ to both points, we may assume $A=(1,0,0,0,0)$
and $B=(x_1, x_2,y, \phi, \eta)\in\hat{L}^+_0$. Note that $x_2$ must be invertible, for otherwise the vectors are not linearly independent in the super-sense, 
namely,  there exist $a, b$ with non-zero bodies so that $aA+bB$ has zero 
body. 

We may apply a transformation $g^s$ of the form \rf{stabil} thus
stabilizing $A$ and mapping $B$ to $g^s\cdot B=(\t x_1, \t x_2,\t y, \t\phi, \t\eta)\in\hat{L}^+_0$, where
\begin{eqnarray}
&&\hskip.4em \t y=y+cx_2+\beta\eta,\quad  {\t x}_2=x_2\nonumber,\\
&&\hskip .4em \t \eta=\beta x_2+\eta,\quad \t\phi=\beta y+c\beta x_2+c\eta+\phi.
\end{eqnarray}
We wish to normalize so that this vector takes the form  
$s(0,1,0, 0, 0)$  
and thus impose the two further conditions $\t \eta=0$ and $\t y=0$. 
This implies
\begin{eqnarray}
\beta=-x_2^{-1}\eta~{\rm and}~c=-\frac{y}{x_2},
\end{eqnarray} 
which gives $s=x_2=\langle A, B\rangle$ as required.
\hfill$\blacksquare$
 
 \begin{Cor}
The moduli space of $OSp(1|2)$-orbits of ordered pairs of vectors in the special light cone lying in distinct rays is given by ${\mathbb R}_+$.
\end{Cor}


Turning now to ordered triples of linearly independent vectors in $\hat{L}^+_0$,  the underlying bosonic triple may provide either a positively or a negatively oriented basis of  $\mathbb{R}^{2,1}$, and we have already  christened the former case a positive triple.
The well-known three-effectiveness of the action of $PSL(2,{\mathbb R})$ on positively oriented (i.e., correctly cyclically ordered) triples of points in the circle at infinity of hyperbolic space fails in our current context of $\mathbb{R}^{2,1|2}$ because there is one additional odd degree of freedom which cannot be fixed. We next
compute the moduli space of orbits of positive triples to be ${\mathbb R}_+^{3|1}$ modulo the fermionic reflection and
postpone the further discussion of this interesting additional parameter to the next section.

\begin{Lem} \label{3pts}Let $\zeta^b \zeta^e \zeta^a$ be a positive triple in the special light cone.
Then 
there is $g\in OSp(1|2)$, which is unique up to composition with the fermionic reflection, and unique even $r,s,t$, which have positive bodies, and odd $\phi$ so that
$$g\cdot \zeta^e=t(1,1,1,\phi,\phi),~g\cdot \zeta^b=r(0,1,0,0,0),~g\cdot \zeta^a=s(1,0,0,0,0).$$
\end{Lem}
\noindent{\bf Proof.} 
First, one can put $\zeta^a$ into the form $(1,0,0,0)$ by means of $OSp(1|2)$ according to Lemma \ref{orb}. Second, $\zeta^b$ can be put 
in the form $(0,f, 0, 0,0)$ using the stabilizer of $(1,0,0,0,0)$ as in the proof of Lemma \ref{2pts}. 
Finally by a diagonal matrix whose entries have positive bodies,  
one can put $\zeta^e$ into the form where its $x_1$- and $x_2$-coordinates agree. Since $\zeta^b \zeta^e \zeta^a$ 
is a positive triple, the $y$-coordinate of $\zeta^e$ after these transformations has positive body, and 
so $\zeta^e$ indeed transforms to a vector with the required form $t(1,1,1,\phi,\phi)$
while the other two transformed vectors evidently also have the desired forms 
$\zeta^b=r(0,1,0,0,0)$, $\zeta^a=s(1,0,0,0,0)$.  The asserted uniqueness up to reflection of the transformation $g\in OSp(1|2)$
thus constructed follows from Corollary \ref{stab}.\hfill ${\blacksquare}$\\
 
Just to be clear, let us note that there are thus precisely two lifts to $\hat L^+_0$ of a positive triple $T$ with abstract coordinates $(r,s,t,\theta)\in{\mathbb R}_+^{3|1}$
given by $g\cdot T$ (where the $y$-coordinate of $\zeta^e$ has positive body in the notation of the lemma) and its image $g^r g \cdot T$
under the fermionic reflection $g^r$.  In any case, the invariants $r,s,t$ have positive body and only the signs of the fermions change.
As a parenthetical point of notation, we mention that the labeling $abe$ is used here rather than the more natural $abc$ in order to accommodate later conventions.

\begin{Cor}\label{modsp3} The moduli space of $OSp(1|2)$-orbits of positive triples in the light cone is given by 
$(r,s,t,\phi)\in{\mathbb R}_+^{3|1}/{\mathbb Z}_2$, where ${\mathbb Z}_2$ acts by fermionic reflection.
\end{Cor}

In fact,  the even invariants $r,s,t$ of a triple, which we shall call {\it normalization coefficients},
are the direct analogues of the reciprocal $h$-lengths from \cite{penner,pb}, and they
can be nicely computed in terms of the 
$\lambda$-lengths
$$a^2=~<\zeta^b,\zeta^e>, ~~b^2=~<\zeta^a,\zeta^e>, ~~e^2=~<\zeta^a,\zeta^b>.$$

\begin{Lem}\label{evenones} The normalization coefficients $r,s,t$ in Lemma \ref{3pts} are given by
\begin{eqnarray}\label{rst}
r=\sqrt{2}~{{ea}\over b},\quad s=\sqrt{2}~{{be}\over a},\quad t=\sqrt{2}~{{ab}\over e}.
\end{eqnarray}
\end{Lem}

We have found that a positive triple of points in the special light cone naturally determines
three even and one odd invariant, namely, the $\lambda$-lengths $a,b,e$ or the normalization coeffcients $r,s,t$ and the fermion label
$\phi$ given in the previous lemma.  It appears that the definition of the odd parameter $\phi$ depends on a choice of member of the positive triple, and we
next discuss how to eliminate this dependence.

In fact given a positive triple $ABC$, there is a canonical cube root  of unity in $OSp(1|2)$ called the {\it $ABC$ prime transformation}
which cyclically permutes its members $(A,B,C)\mapsto (B,C,A)$,
and we shall compute this transformation on the coordinates of Lemma \ref{3pts} explicitly and denote it
$(r,s,t, \phi)\mapsto (r',s',t',\phi')$ where the $\lambda$-lengths are of course also cyclically permuted in
the natural way.  The choice-free way to express the additional odd degree of freedom
is to average our ansatz $\phi$ over this prime transformation and define the  {\it Manin~invariant}\footnote{In fact, Manin \cite{manin} introduced the odd ``pseudo-invariant'' $\pm\theta$ of a triple of points in ${\mathbb R}^{1|1}$
capturing the basic non-transitivity of the $OSp(1|2)$ action there,
and by making choices (of spin structure among other conventions) one can \cite{zababaher}  lift this
to a signed expression $\theta$, namely, our Manin invariant.} or simply the {\it $\mu$-invariant}
to be
\begin{eqnarray}\label{manavg}
&\theta ={1\over 3}(\phi +\phi'+\phi'').
\end{eqnarray}
In practice, $\mu$-invariants and $\lambda$-lengths give a parametrization of decorated super-Teichm\"uller space as we shall see.

\begin{figure}[hbt] 
\centering           
\includegraphics[width=0.3\textwidth]{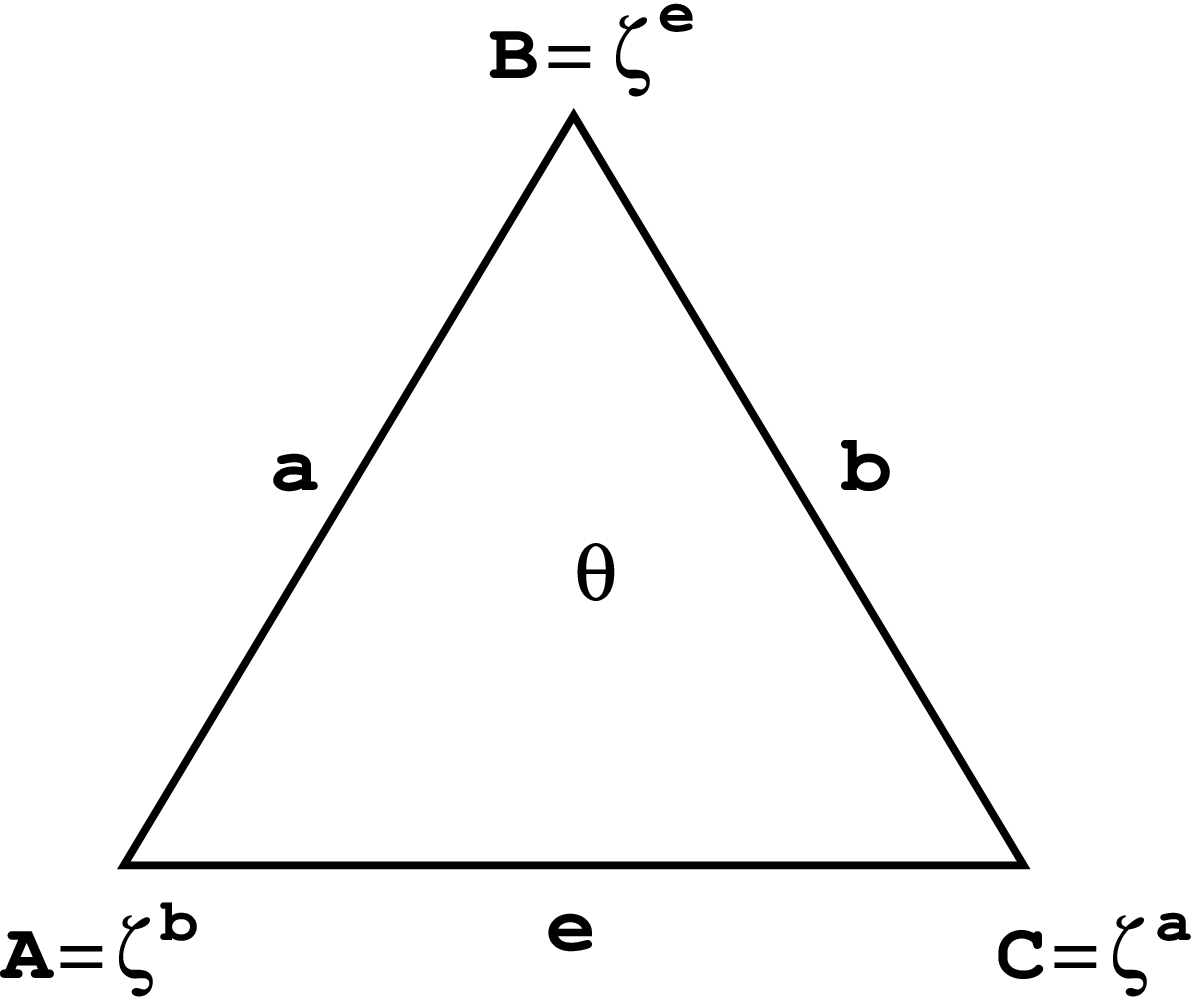}
\caption{Parametrization of a positive triple in the special light cone by three $\lambda$-lengths $a,b,e$ on edges and one $\mu$-invariant $\theta$.}\label{3nota}                
\end{figure}

To complete our definition of the $\mu$-invariant, it therefore
remains to compute the prime transformation.  To this end, we alter notation
slightly setting $A=\zeta^b$, $B=\zeta^e$, $C=\zeta^a$ as illustrated in
Figure \ref{3nota}.

Given three points 
\begin{eqnarray}
&&A=r(0,1,0,0, 0)\nonumber,\\
&&B=t(1,1,1, \phi, \phi),\nonumber\\
&&C=s(1,0,0,0,0)
\end{eqnarray}
in the special light cone $\hat{L}^+_0$,
where $r,s,t$ have positive body, we claim there exists a group element $g^{A',B',C'}_{A,B,C}$ from $OSp(1|2)$    
respectively transforming the points $A,B,C$ into  
\begin{eqnarray}
&&A'=t'(1,1,1, \phi', \phi'),\nonumber\\
&&B'=s'(1,0,0,0,0),\nonumber\\
&&C'=r'(0,1,0,0, 0),
\end{eqnarray}
and in fact, there are two possible ways to construct such a 
transformation (again unique up to fermionic reflection).
However, there is a unique transformation whose third power
\begin{eqnarray}\label{gid}
g^{A''',B''',C'''}_{A'',B'',C''}
g^{A'',B'',C''}_{A',B',C'}g^{A',B',C'}_{A,B,C}=1
\end{eqnarray}
is equal to identity, and it turns out that its action on $\phi$ is rather simple. 
\begin{Prop} \label{prime} 
The group element $g^{A',B',C'}_{A,B,C}$ acts so that $\phi''=\phi'=\phi$, and this is the unique group element $g^{A',B',C'}_{A,B,C}\in OSp(1|2)$ satisfying the identity \rf{gid}.  
\end{Prop}

\noindent{\bf Proof.} We first use  
the group element $g_{\phi\psi}^t$ in \rf{gtpp} in order to transform $B$ into $(1,0,0,0,0)$. The effects of this transformation on the other points produces
\begin{eqnarray}
&&g^t_{\phi\phi}: C\mapsto st(0, 1,0,0, 0),\nonumber\\
&&g^t_{\phi\phi}: A\mapsto r(t^{-1}, t, 1, -\frac{\phi}{\sqrt{t}}, -\sqrt{t}\phi).
\end{eqnarray}
Finally, acting by the diagonal matrix $g^d=\biggl(\begin{smallmatrix}
\sqrt{t}&0&0\\
0&\sqrt{t^{-1}}&0\\
0&0&1\\
\end{smallmatrix}\biggr)$
we obtain the desired result by composing with the fermionic reflection $g^f$
\begin{eqnarray}
g^{A',B',C'}_{A,B,C}=g^t_{\phi,\phi}g^dg^f=\biggl(\begin{smallmatrix}
\hskip1ex0&\hskip1ex1&\hskip1ex0\\
-1&-1&-\phi\\
\hskip1ex0&-\phi&\hskip1ex1\\
\end{smallmatrix}\biggr)
\end{eqnarray}
so that $t'=r$, $r'=s$, $s'=t$. $\blacksquare$ \\

\begin{Cor} The $\mu$-invariant $\theta\equiv \phi$ of a positive triple $ABC$ is invariant under the $ABC$ prime transformation.
\end{Cor}

As follows directly from this plus the discussion of the previous section, we have

\begin{Thm} The moduli space of $OSp(1|2)$-orbits of positive triples in the light cone is given by equivalence classes under fermionic reflection of
three even $\lambda$-lengths $a,b,e$ with positive bodies plus the $\mu$-invariant $\theta$ defined up to fermionic reflection. 
\end{Thm}

\section{Basic calculation and Ptolemy transformations}\label{sec:bcpt}
This section is dedicated to a computation called the ``basic calculation''
giving a para\-met\-rization of the moduli space of $OSp(1|2)$-orbits of
four-tuples $ABCD$ in the special light cone $\hat{L}_0^+$ comprised of two positive triples $CBA$ and $DCA$ of points
in terms of five $\lambda$-lengths and two $\mu$-invariants.  In effect, one
positive triple $CBA$ is put into the canonical position of Lemma \ref{3pts} using certain of the
putative parameters, and the remaining ones are then used to completely and uniquely determine 
the fourth point $D=(x_1,x_2,-y,\rho,\lambda)$, where $y$ has non-negative body; it is precisely here that compatibility of signs of adjacent $\mu$-invariants arises.  This basic calculation is the critical ingredient both for our global coordinates
and for the Ptolemy transformation just as in \cite{penner,pb}.
In fact provided only that the triple $CBA$ is in standard position, the standard position of $DCA$ is easily calculated in terms of the 
coordinates of its vertices as follows.

\begin{Lem}\label{hats} The triple
\begin{eqnarray}
&A=r(0,1,0,0,0), ~C=s(1,0,0,0,0), ~D=(x_1,x_2,-y,\rho,\lambda), ~{\rm for}~y\geq 0,\nonumber
\end{eqnarray}
can be transformed by an element of $OSp(1|2)$, which is itself uniquely determined up to fermionic reflection, into the triple
\begin{eqnarray}
&\hat A=\hat s(1,0,0,0,0), ~~\hat C=\hat r(0,1,0,0,0)), ~~\hat D=\hat{t}(1,1,1,\sigma,\sigma),\nonumber
\end{eqnarray}
where 
\begin{eqnarray}
&&\hat s=\sqrt{\frac{x_1}{x_2}} r,\quad \hat r=\sqrt{\frac{x_2}{x_1}}s,\quad 
\hat t =\sqrt{x_1x_2},\nonumber\\ 
&&\sigma=-\sqrt[4]{\frac{x_1}{x_2}}\frac{\lambda}{\sqrt{x_1x_2}}=\sqrt[4]{\frac{x_2}{x_1}}\frac{\rho}{\sqrt{x_1x_2}},
\end{eqnarray}
\end{Lem}

\noindent{\bf Proof.}  Act first by $\left( \begin{smallmatrix}
\hskip1ex0 & 1 & 0 \\
-1 & 0 & 0 \\
\hskip1ex0 & 0 & 1 \end{smallmatrix} \right )$ to interchange the roles of $A$ and $C$ while mapping $D$ to a point with non-negative $y$-coordinate and then act by $\left( \begin{smallmatrix}
\sqrt[4]{{x_1}/ {x_2}} & 0 & 0 \\
0 & \sqrt[4]{{x_2}/ {x_1}} & 0 \\
0 & 0 & 1 \end{smallmatrix} \right)$ to realize equality of the first two coordinates.  
The composition of these two elements of $OSp(1|2)$ achieves the desired result in accordance with Lemma \ref{3pts}.
\hfill $\blacksquare$\\

\noindent Note that the normalization factors admit the uniform expressions
\begin{eqnarray}\label{2rst}
&&t=\sqrt{2}\frac{a b}{e},\quad
r=\sqrt{2}~{{ea}\over b},\quad 
s=\sqrt{2}~{{be}\over a},\nonumber\\
&&\hat t=\sqrt{2}\frac{cd}{e},\quad
\hat r=\sqrt{2}~{{ec}\over d},\quad 
\hat s=\sqrt{2}~{{de}\over c}.
\end{eqnarray}
This transformation in $OSp(1|2)$ which maps $DCA$ to standard position if $BAC$ is in standard position is called the {\it switch transformation} of the quadruple $ABCD$ though in fact it depends only upon $DCA$ according to the formula for it given above.

\begin{figure}[hbt] 
\centering           
\includegraphics[width=0.6\textwidth]{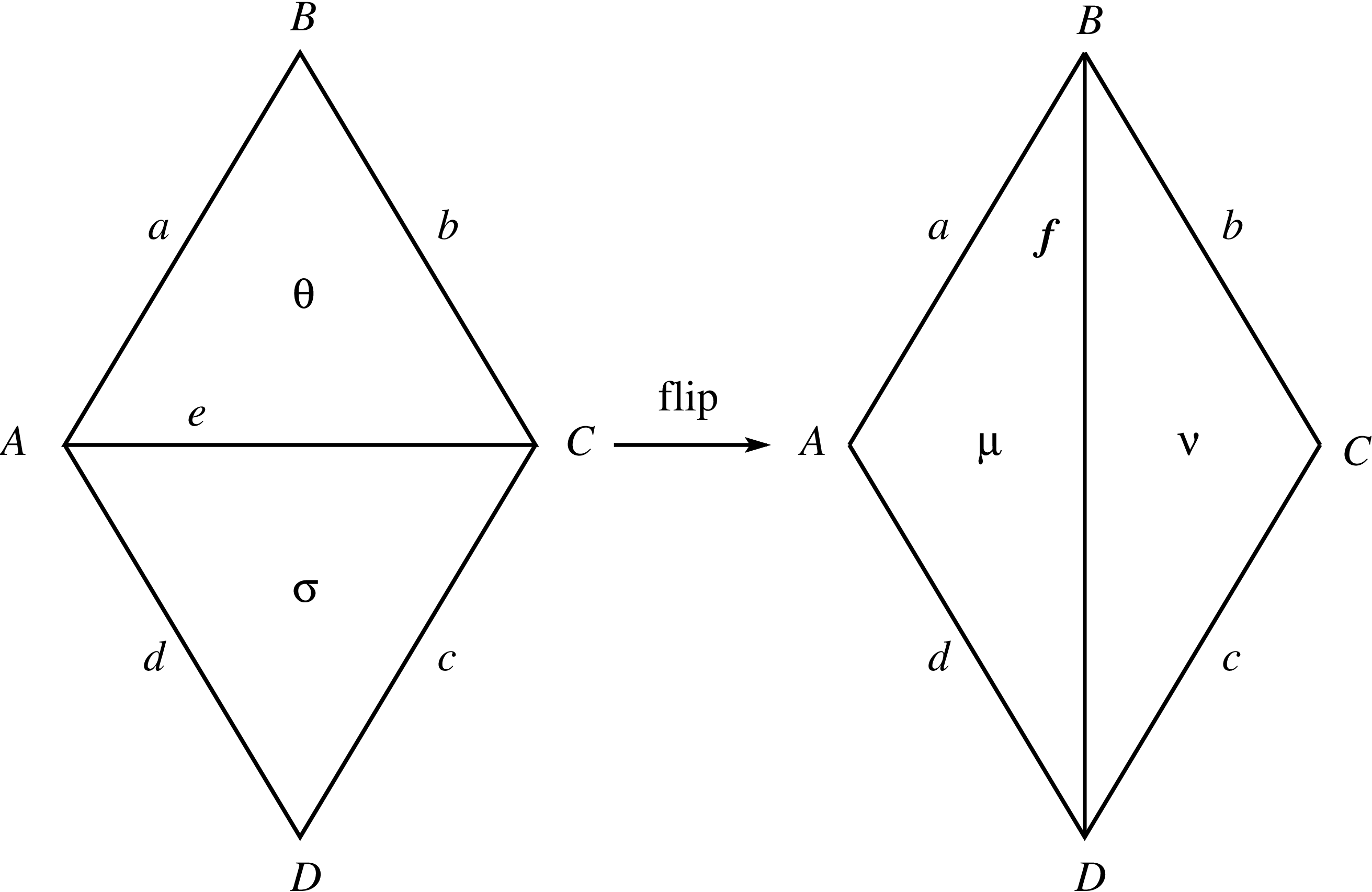}
\caption{Standard notation for $\lambda$-lengths and $\mu$-invariants near an edge before (left) and after (right) a flip}  \label{figint}              
\end{figure}


Now turning to the basic calculation itself, let $A=\zeta^b, B=\zeta^e, C=\zeta^a$ be points in the light cone as in Lemma \ref{3pts} with
$\lambda$-lengths $a,b,e$ and
$\mu$-invariant $\theta$. Up to fermionic reflection, the $OSp(1|2)$-orbit of the positive triple $CBA$ is then uniquely determined
by these parameters according to Lemma \ref{3pts}.
The basic calculation aims to compute $D=(x_1,x_2, -y, \rho, \lambda)\in\hat{L}^+_0$, where $y$ has non-negative body,
from the data
\begin{eqnarray}\label{4point}
&&\sigma=-\sqrt[4]{\frac{x_1}{x_2}}\frac{\lambda}{\sqrt{x_1x_2}}=\sqrt[4]{\frac{x_2}{x_1}}\frac{\rho}{\sqrt{x_1x_2}},
\quad\langle A,D\rangle=d^2, \quad \langle D, C\rangle=c^2
\end{eqnarray}
depicted in Figure \ref{figint}, and the following arises directly from equations \rf{rst} and  \rf{4point}.

\begin{Prop}\label{bascalc}
We have the expressions
\begin{eqnarray}
&&x_1=\frac{\hat t \hat s}{r}=\sqrt{2}\frac{d^2b}{ea}, \quad x_2=\frac{\hat r\hat t}{s}=\sqrt{2}\frac{c^2a}{eb},\nonumber\\ 
&&y=\sqrt{x_1x_2}=\sqrt{2}\frac{cd}{e},\nonumber\\
&&\lambda=-\sqrt{2}\frac{cd}{e}\sqrt{\frac{ac}{db}}~\sigma, \quad \rho=\sqrt{2}
\frac{cd}{e}\sqrt{\frac{db}{ac}}~\sigma
\end{eqnarray}
solving the basic calculation.
\end{Prop}

\noindent Introducing the cross-ratio $Z=\frac{ac}{db}$, this expression of the basic calculation can be concisely written
\begin{eqnarray}
&&x_1=\sqrt{2}\frac{cd}{e}Z^{-1}, 
\quad x_2=\sqrt{2}\frac{cd}{e}Z,\quad y=\sqrt{2}\frac{cd}{e},\nonumber\\
&&\lambda=-\sqrt{2}\frac{cd}{e}\sqrt{Z}~\sigma, \quad \rho=\sqrt{2}
\frac{cd}{e}\sqrt{Z^{-1}}~\sigma.
\end{eqnarray}

The utility of this version of the basic calculation for computing the Ptolemy transformation is already evident from the following proposition.

\begin{Prop}\label{pt*} In the notation above with
$f^2=\langle A,D\rangle$,  we have
\begin{eqnarray}
&&ef=(ac+bd)\Big(1+
\frac{\sigma\theta}{{Z}^{1\over 2}+{Z}^{-{1\over 2}}}\Big).
\end{eqnarray} 
\end{Prop}
\noindent{\bf Proof.} Simply expand the expression
\begin{eqnarray}
&&f^2=\frac{t}{2}(x_1+x_2)+ty+ t\rho\theta+t\theta\lambda\nonumber\\
&&\hskip .2in =\frac{d^2b^2}{e^2}+\frac{a^2c^2}{e^2}+
2\frac{abcd}{e^2}+2\frac{abcd}{e^2}\Big(\sqrt{\frac{db}{ac}}+
\sqrt{\frac{ac}{db}}\Big)\sigma\theta
\end{eqnarray} 
in $\lambda$-lengths.
This is equivalent to 
\begin{eqnarray}
e^2f^2=(ac+bd)^2+2abcd({{Z}^{1\over 2}+{Z}^{-{1\over 2}}})\sigma\theta,
\end{eqnarray} 
and the desired formula for $f$ follows upon taking square roots.\hfill $\blacksquare$\\

Now let us reproduce the same four points $A, B, C, D\in\hat{L}^+_0$ using the parameters of adjacent positive triples $CBD$ and $ADB$. To do so, let us move the positive triple $CBD$ to standard position
\begin{eqnarray}\label{tildas}
&&B=
\t s (1,0,0,0,0), 
\quad C=\t r (0,1,0,0, 0),\nonumber
\\
&&D=\t t(1,1,1, \t \phi, \t \phi).
\end{eqnarray}
To this end, we first apply the stabilizer of $s(1,0,0,0,0)$, namely, the matrix 
\begin{eqnarray}\label{stabpt}
g=\left( \begin{array}{ccc}
1 & 0 & 0 \\
c & 1 & \zeta \\
\zeta & 0 & 1 \end{array} \right),
\end{eqnarray}
as in Corollary \ref{stab},
in order to transform $t(1,1,1, \theta, \theta)\to \bar r(0,1,0,0, 0)$, where the corresponding $c$ and $\zeta$ are given by
$
c=-1~{\rm and}~ \zeta=-\theta 
$
yielding $\bar{r}=t$.
Applying this transformation to 
$D=(x_1, x_2, -y, \rho, \lambda)\in\hat{L}^+_0$, we obtain a new vector $(\hat{x}_1, \hat{x}_2, -\hat{y}, \hat{\rho}, \hat{\lambda})\in\hat{L}^+_0$, 
where
\begin{eqnarray}
&&\hat{x}_1=x_1+x_2+2y-2\theta(\rho-\lambda),\nonumber\\
&&\hat{x}_2=x_2,\nonumber\\
&&\hat{\rho}=\theta y+\theta x_2-\lambda+\rho,\nonumber\\
&&\hat{\lambda}=-\theta x_2+\lambda.
\end{eqnarray}

Applying $\bigl (\begin{smallmatrix}\hskip1ex 0&1&0\\-1&0&0\\\hskip1ex0&0&1\end{smallmatrix}\bigr)$ followed by the diagonal matrix 
$\biggl (\begin{smallmatrix}\sqrt[4]{\frac{\hat{x}_1}{\hat{x}_2}}&0&0\cr 0&\sqrt[4]{\frac{\hat{x}_2}{\hat{x}_1}}&0\cr 0&0&1\end{smallmatrix}\biggr)$
to each of $B,C,D$ yields:
\begin{eqnarray}\label{bcd}
&&\t t=\sqrt{\hat{x}_1\hat{x}_2},\quad
\tilde r=t\sqrt{\frac{\hat{x}_1}{\hat{x}_2}},\quad \tilde s=s \sqrt{\frac{\hat{x}_2}{\hat{x}_1}}\nonumber\\
&&\tilde\phi=-\frac{\hat{\lambda}}{\sqrt{\hat{x}_1\hat{x}_2}}
\Big(\sqrt[4]{\frac{\hat{x}_1}{\hat{x}_2}}\Big)=\frac{\hat{\rho}}{\sqrt{\hat{x}_1\hat{x}_2}}
\Big(\sqrt[4]{\frac{\hat{x}_2}{\hat{x}_1}}\Big)
\end{eqnarray}
 for the parameters in equation \rf{tildas}.
Thus, the $\mu$-invariant of the positive triple $BDC$ is given by 
$\t\phi$ as this is invariant under the $BDC$ prime transformation.

\begin{Thm}\label{fl1}
The $\mu$-invariants for the positive triples $BDC$  and $DBA$
depend only on the $\mu$-invariants $\theta, \sigma$ and the cross-ratio $Z$ and are given by
\begin{eqnarray}\label{one}
\nu=\frac{\theta\sqrt{Z}+\sigma}{\sqrt{1+Z}},
\end{eqnarray} 
\begin{eqnarray}\label{two}
\mu=\frac{\sigma\sqrt{Z}-\theta}{\sqrt{1+Z}}.
\end{eqnarray}
\end{Thm}
\noindent{\bf Proof.} Equation \rf{one} follows from the direct calculation
\begin{eqnarray}\label{three}
\nu=-\frac{\hat{\lambda}}{\sqrt{\hat{x}_1\hat{x}_2}}
\Big(\sqrt[4]{\frac{\hat{x}_1}{\hat{x}_2}}\Big)=\frac{\theta x_2-\lambda}
{Z^{\frac{3}{4}}(\sqrt{2}\frac{cd}{e})\sqrt{\sqrt{Z}+\sqrt{Z}^{-1}}}=
 \frac{\theta\sqrt{Z}+\sigma}{\sqrt{1+Z}}.
\end{eqnarray}
In order to write $\mu$ in terms of the same function $f^{BDC}$, we must
transform $D, A, B\in\hat{L}^+_0$ as follows
\begin{eqnarray}\label{tcheck}\nonumber
&&D\mapsto \check{t}(1,1,1,\check{\phi},\check{\phi}),\quad 
A\mapsto \check{s}(1,0,0,0,0)\nonumber,\quad
B\mapsto(\check{x}_1, \check{x}_2, -\check{y},\check{\rho},\check{\lambda}),
\end{eqnarray}
and this can be achieved by first applying the matrices given in Lemma \ref{hats} so that  
\begin{eqnarray}
&&\check{\rho}=-t\sqrt[4]{\frac{x_1}{x_2}}\theta, \quad \check \lambda=
t\sqrt[4]{\frac{x_2}{x_1}}\theta,\nonumber\\
&&\check\phi=-\sqrt[4]{\frac{x_1}{x_2}}\frac{\lambda}{\sqrt{x_1x_2}}=\sqrt[4]{\frac{x_2}{x_1}}\frac{\rho}{\sqrt{x_1x_2}}.\nonumber
\end{eqnarray}
It follows that 
\begin{eqnarray}
\mu=f^{BDC}(\check{Z}, \check{\theta}, \check{\sigma}),
\end{eqnarray}
where the checked arguments of $\mu$ are related to checked variables from \rf{tcheck} as before and
$
\check{Z}=Z, \check\theta=\sigma, \check\sigma=-\theta.
$
It follows that $\mu=f^{DBA}(Z,\theta,\sigma)=f^{BDC}(Z,\sigma,-\theta)$
as required.
\hfill $\blacksquare$

\medskip



As a direct corollary to the proof, we have

\begin{Cor}\label{switch}
Consider the odd Ptolemy transformation on ordered pairs
$(\theta, \sigma)\mapsto (\nu,\mu)$ together
with the corresponding action on $\lambda$-lengths
and apply it twice
to the quadrilateral $ABCD$, with $BAC$ in 
standard position.  
Then the effect is the {switch transformation} of $ABCD$
described in Lemma \ref{hats}. 
Specifically after this transformation, $DCA$ is in standard position with $D\mapsto \hat{t}(1,1,1,\sigma,\sigma)$ while the image of $B$ is determined not by $\theta$ but rather by $-\theta$.
\end{Cor}

\section{Spin surfaces and orientations on fatgraphs}
In this section, we relate the collection of spin structures on a punctured surface
to orientations on any trivalent fatgraph spine of the surface.   
To begin,
we recall results of Cimasoni-Reshetikhin from \cite{cr1,cr2}.

The {\it boundary}  of a one-dimensional CW complex $\mathcal{G}$ is its set $\partial{\mathcal G}$ of vertices of valence one.
$\mathcal{G}$ is a {\it surface graph with boundary} for some compact oriented surface $\Sigma$ with boundary $\p \Sigma$ if $\mathcal{G}$ is embedded in $\Sigma$ with 
$\mathcal{G}\cap \p\Sigma=\p\mathcal{G}$ so that $\b{\mathcal{G}}=\mathcal{G}\cup \p \Sigma$ is the 1-skeleton of a cellular decomposition of $\Sigma$. 
A {\it dimer configuration} on $\mathcal{G}$ is a choice of certain of its edges called {\it dimers} such that each vertex not in $\p{\mathcal G}$ has exactly one incident dimer while vertices in
$\p{\mathcal G}$ may or may not have an incident dimer, the specification of which is 
regarded as an {\sl a priori} boundary condition on the dimer configuration.

Given an orientation $K$ on the edges of $\bar{\mathcal G}$ and a closed oriented edge curve $C$ in it, we denote by $n^K_C$ the number of edges counted with multiplicity where the orientation of $C$ disagrees with that of $K$.
A {\it Kasteleyn orientation} on $\b{\mathcal{G}}$ is an orientation $K$ on the edges of 
the 1-skeleton $\bar{\mathcal G}$ so that 
$n^K_{\p f}=1({\rm mod}~2)$ for each face $f$ of $\Sigma$.   
Define an equivalence relation $K_1\sim K_2$ between two Kasteleyn orientations $K_1,K_2$ on ${\mathcal G}$ generated by altering the orientation on every edge incident on some fixed vertex, which is called a {\it Kastelyn reflection}, and let ${\mathcal K}({\mathcal G})$ denote the set of equivalence classes.
There is furthermore a cochain $\theta_{K_1,K_2}\in C^1(\Sigma; \mathbb{Z}_2)$ defined by setting $\theta_{K_1,K_2}(e)=1$ if and only if $K_1$ and $K_2$ disagree on $e$, for $e$ an edge of ${\mathcal G}$.

\begin{Thm} [Corollary 1 of \cite{cr1}]\label{theta}
In fact, $\theta_{K_1, K_2}\in Z^1(\Sigma; \mathbb{Z}_2)$ is a cocycle which is a coboundary if and only if  $K_1\sim K_2$.  Furthermore, the set of equivalence classes of Kasteleyn orientations is an affine $H^1(\Sigma; \mathbb{Z}_2)$-space.
\end{Thm}

Seminal work of Dennis Johnson \cite{johnson} mentioned before identifies as
affine $H^1(\Sigma;\mathbb{Z}_2)$-spaces the collection $\Omega(F)$ of all spin structures on a surface $\Sigma$ with the collection ${\mathcal Q}(\Sigma)$ of all quadratic forms
$q: H_1(\Sigma; \mathbb{Z}_2)\to \mathbb{Z}_2$ satisfying
$q(a+b)=q(a)+q(b)+a\cdot b$ with respect to the homology intersection pairing $a\cdot b$ of $a,b\in H_1(\Sigma; {\mathbb Z}_2)$.
In fact, Kasteleyn orientations on a surface graph for $\Sigma$ and quadratic functions on $H_1(\Sigma;{\mathbb Z}_2)$ are also isomorphic as affine $H^1(\Sigma; \mathbb{Z}_2)$-spaces, and the correspondence is given via an explicit construction relative to a fixed dimer $D$ as follows.

\begin{Thm}[Theorem 2.2 of \cite{cr2}]\label{cr}
Fix a dimer configuration $D$ on a surface graph with boundary $\mathcal{G}$ for  the surface $\Sigma$ and let 
$\alpha\in H_1(\Sigma; \mathbb{Z}_2)$ be represented by oriented closed curves $C_1,\dots, C_m\in \bar{\mathcal{G}}$. If $K$ is a Kasteleyn orientation on $\mathcal{G}$, then the function $q^K_D: H_1(\Sigma; \mathbb{Z}_2)\to \mathbb{Z}_2$ given by 
\begin{eqnarray}\label{q}
q^K_{D}(\alpha)=\sum_{i<j}C_i\cdot C_j+\sum^m_{n=1}(1+n^K_{C_i}+\ell^D_{C_i})\quad  ({\rm mod} ~2)
\end{eqnarray}  
is a well-defined quadratic form,
where $\ell^D_C$ is the number of edges of $D$ sticking out to the left of $C$,
and
$n^K_C$ as before is the number of edges counted with multiplicity where the orientation of $C$ disagrees with that of $K$.  Moreover, for each fixed dimer $D$, this establishes
an isomorphism ${\mathcal K}({\mathcal G})\approx{\mathcal Q}(\Sigma)$ as affine $H^1(\Sigma; \mathbb{Z}_2)$-spaces.\end{Thm}

\noindent Together with \cite{johnson}, this establishes an isomorphism of affine
$H^1(\Sigma; \mathbb{Z}_2)$-spaces between the collection $\Omega(\Sigma)$ of spin structures  on $\Sigma$ and ${\mathcal K}({\mathcal G})$ for any surface graph ${\mathcal G}$ with boundary for $\Sigma$ with respect to a fixed dimer configuration.

\begin{figure}[hbt] 
\centering           
\includegraphics[width=0.4\textwidth]{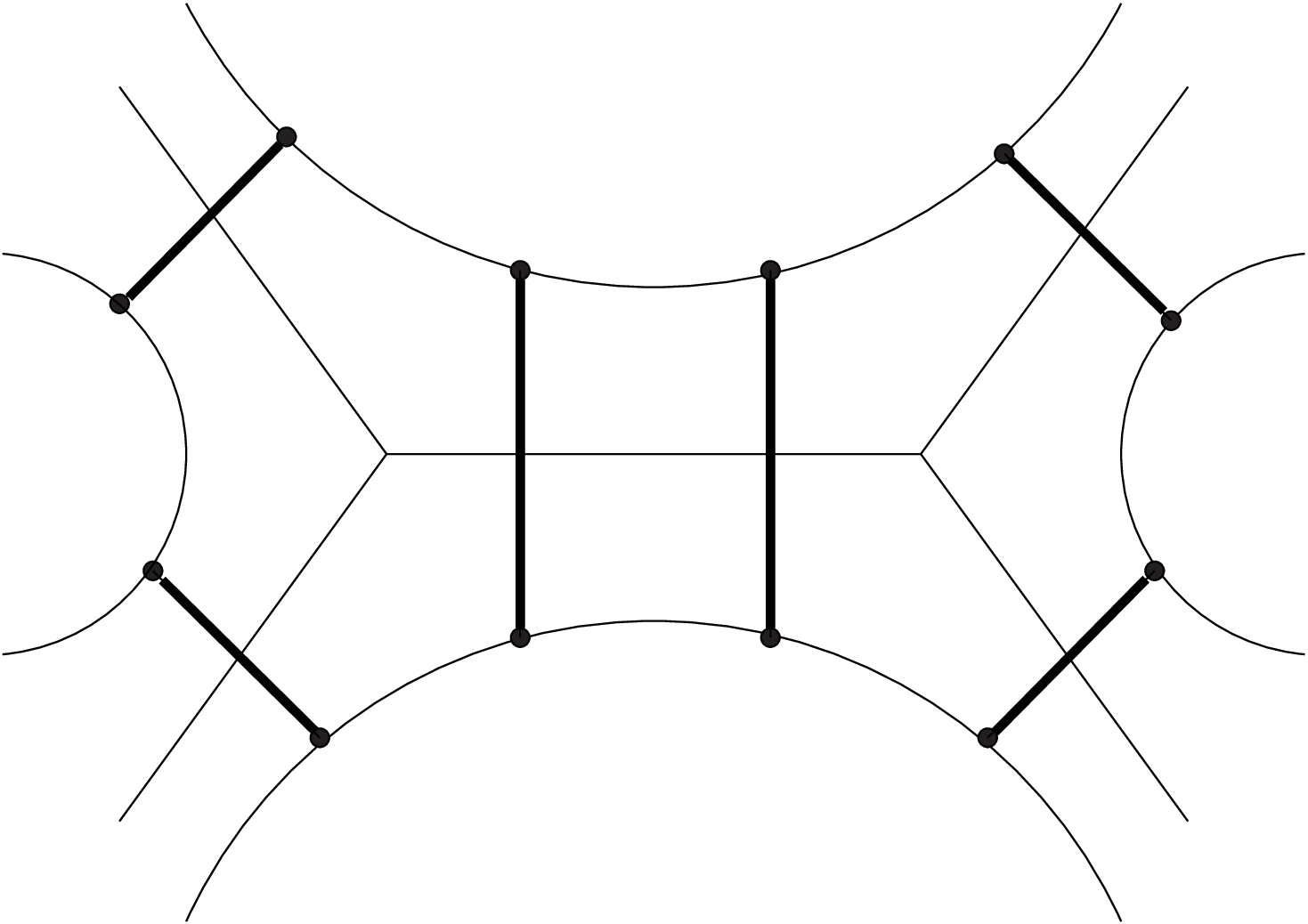}\quad\quad  
\includegraphics[width=0.4\textwidth]{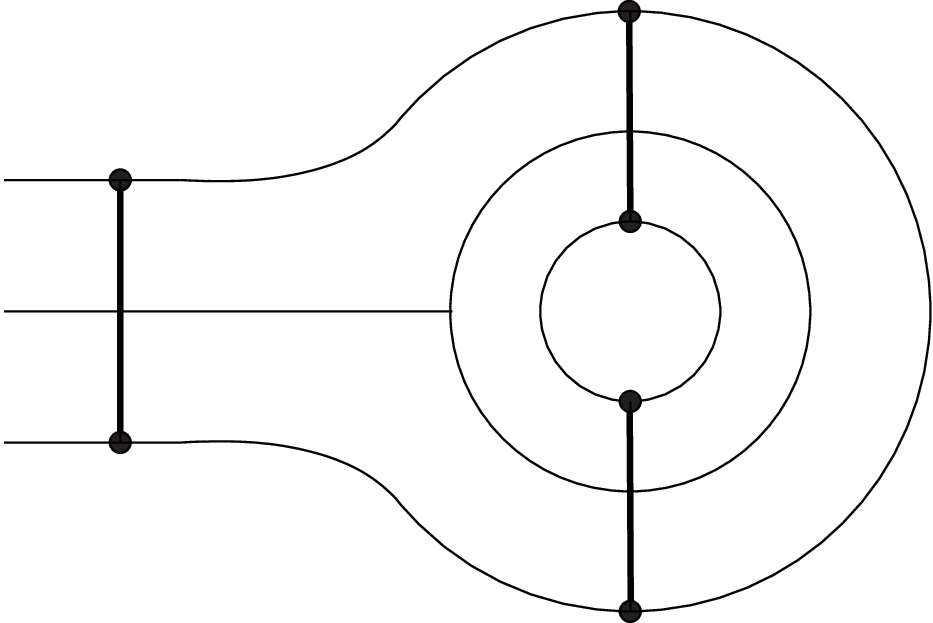}
\caption{Cellular decomposition via fatgraphs}  \label{cell}               
\end{figure}

Now, given a trivalent fatgraph spine $\tau$ for $F=F_g^s$, we shall build
an appropriate surface graph with boundary for a surface embedded in $F$ as a deformation retract.   Construct a CW decomposition of this compact so-called skinny surface $\Sigma=\Sigma(\tau)\subset F$ with boundary taking one hexagon $H_v$ for each vertex $v$ and one rectangle $R_e$ for each edge $e$ of $\tau$ glued together in the natural way in $F$ as in Figure \ref{cell}.
There is a canonical surface graph ${\mathcal{G}}={\mathcal G}_\tau$ for $\Sigma$ comprised of the common boundaries
of these hexagonal and rectangular regions as also illustrated in the figure by bold line segments, two such
segments for each edge of $\tau$.  There is also a canonical dimer ${\mathcal D}=\mathcal{G}$ given by exactly this same set of segments, so each vertex in $\p{\mathcal G}$  has an incident dimer.

\begin{figure}[hbt] 
\centering           
\includegraphics[width=0.8\textwidth]{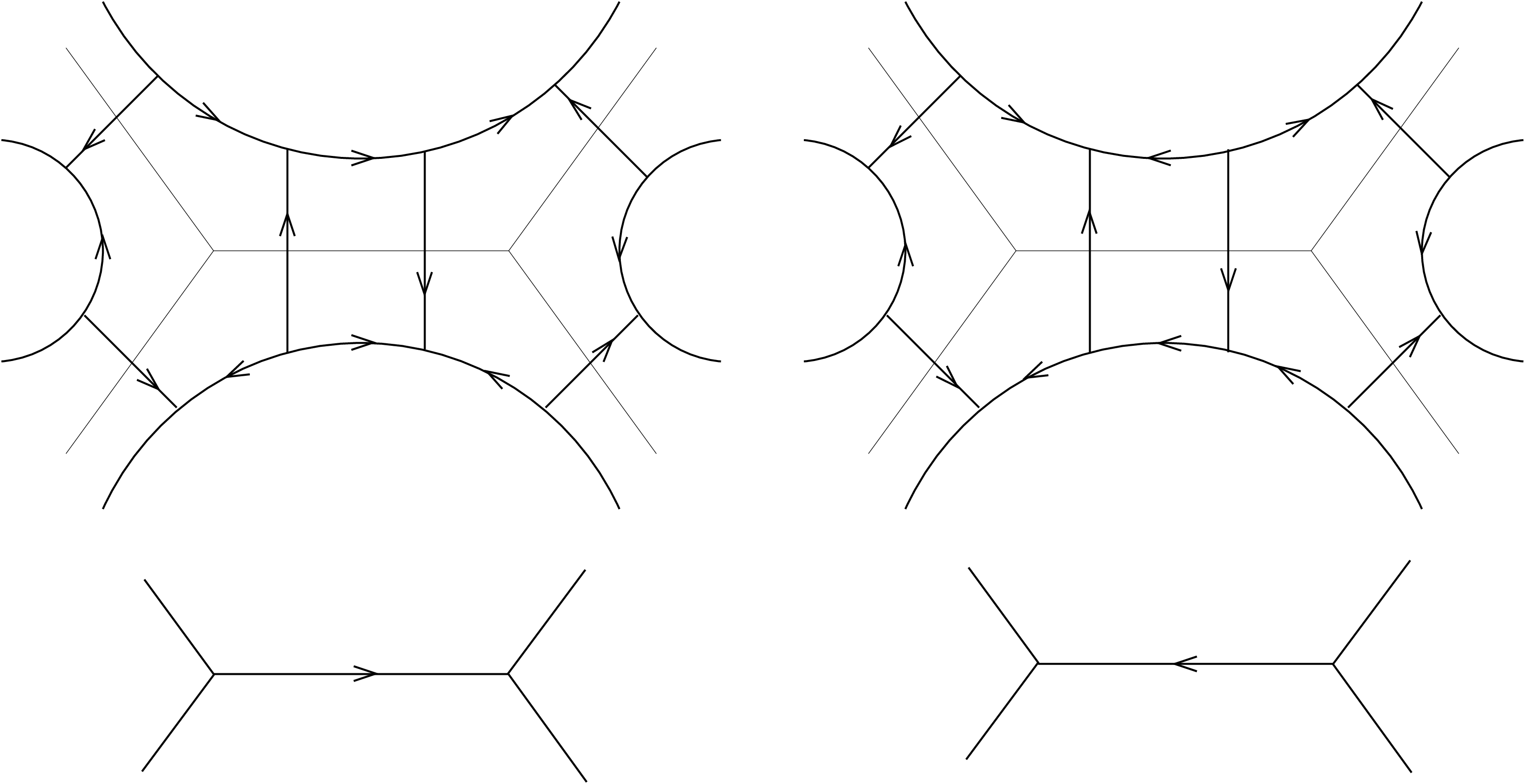}
\caption{Extension of special Kastelyn orientations on hexagons to special ones on $\b{\mathcal G}$ on the top and orientation on fatgraph $\tau$ from special Kasteleyn orientations on $\bar{\mathcal{G}}$ on the bottom} \label{orient}               
\end{figure}

A hexagon has the two special Kastelyn orientations that are either outgoing or incoming at each vertex, and these are related by reversal of orientation of each edge.
Furthermore identifying $H_v$ with this abstract hexagon, for some vertex $v$ of $\tau$, there is a unique such orientation which agrees with the one induced
from the counter-clockwise orientation of ${\mathcal G}$ with $H_v$ on the left
as illustrated on the top of Figure \ref{orient}.
Thus, any Kastelyn orientation on $\b{\mathcal G}$ can be modified by Kastelyn reflections to agree with this one on $H_v$
for each vertex $v$ of $\tau$, i.e., the special Kastelyn orientations saturate the equivalence classes.



Suppose that $K$ is such a special Kastelyn orientation on $\b{\mathcal G}$ that thus agrees with
the special hexagonal ones near each vertex of $\tau$ oriented as already discussed.  There are exactly two methods that $K$ can extend as a Kastelyn orientation to a rectangle $R_e$  as illustrated on the top in Figure \ref{orient}, and these are naturally in bijective correspondence with orientations on the edges of $\tau$ themselves insofar as they are parallel as illustrated on the bottom of Figure \ref{orient}.

Conversely, suppose that $\omega$ is an orientation on $\tau$ and let $K_\omega$ denote the associated special Kastelyn orientation on $\b{\mathcal G}$ defined to be parallel on each $R_e$ and special on each $H_v$ as before.   
In the same spirit as Kastelyn reflection, we define a {\it fatgraph reflection} at a vertex $v$ of $\tau$ to reverse the orientations of $\omega$ on every edge of $\tau$ incident on $v$, consider the equivalence relation $\omega_1\sim\omega_2$ thus generated on the set of all orientations on $\tau$ and let ${\mathcal O}(\tau)$ denote the set of all equivalence classes.  One easily checks that the fatgraph reflection at $v$ is given precisely by the composition of Kastelyn reflections at the six vertices of $H_v$, so $\omega_1\sim\omega_2$ implies $K_{\omega_1}\sim K_{\omega _2}$.
Conversely, in order that the orientation on an edge of some $H_v$ is invariant, we must
perform Kastelyn reflection either at both or neither of its endpoints, and it follows that
$K_{\omega_1}\sim K_{\omega _2}$ also implies $\omega_1\sim\omega_2$.  Furthermore, given
orientations $\omega_1,\omega_2$ on $\tau$, there is the analogous cochain
$\theta_{\omega_1,\omega_2}\in C^1(\Sigma;{\mathbb Z}_2)$ taking a non-zero value on an edge of $\tau$ if and only if $\omega_1$ and $\omega_2$ disagree on the edge, and we clearly have $\theta_{\omega_1,\omega_2}\equiv\theta_{K_{\omega_1},K_{\omega_2}}$.
Summarizing, we have:


\begin{Prop}\label{equiva}
For each orientation $\omega$ on the edges of a fatgraph $\tau$, there exists a unique special Kasteleyn orientation $K_{\omega}$ on $\bar{\mathcal{G}}_{\tau}$, and this induces an isomorphism 
$ {\mathcal K}({\mathcal G}_\tau)\approx\mathcal{O}({\tau})$ of equivalence classes under reflection as affine $H^1(F^s_g; \mathbb{Z}_2)$-spaces.
\end{Prop}

As follows directly from Theorem \ref{cr} and Proposition \ref{equiva}, we have

\begin{Thm}\label{qando}
Let $\tau$ be a fatgraph  spine in the surface $F=F_g^s$ with corresponding surface graph ${\mathcal G}_\tau$ for $\Sigma\subset F$.  Then the formula \rf{q}
with the canonical dimer configuration ${\mathcal D}$
establishes an isomorphism ${\mathcal O}(\tau)\approx {\mathcal Q}(\Sigma)$
as affine $H^1(F;{\mathbb Z}_2)$-spaces, and indeed also
$ {\mathcal Q}(\Sigma)\approx {\mathcal Q}(F)$ since $\Sigma\subset F$ is a homotopy equivalence.
\end{Thm}

It remains for us here only to compute the effect that flipping a fatgraph edge has on 
an orientation class:

\begin{Lem}
Suppose that $\tau_1$ is a trivalent fatgraph spine for $F$ and that $\tau_2$
arises by flipping an edge of $\tau_1$.
There is a unique bijection ${\mathcal O}(\tau_1)\to{\mathcal O}(\tau_2)$ covering the identity map of ${\mathcal Q}(F)$, and it is described by Figure \ref{flipper}.
\end{Lem}

\smallskip

\begin{figure}[hbt] 
\centering           
\includegraphics[width=0.4\textwidth]{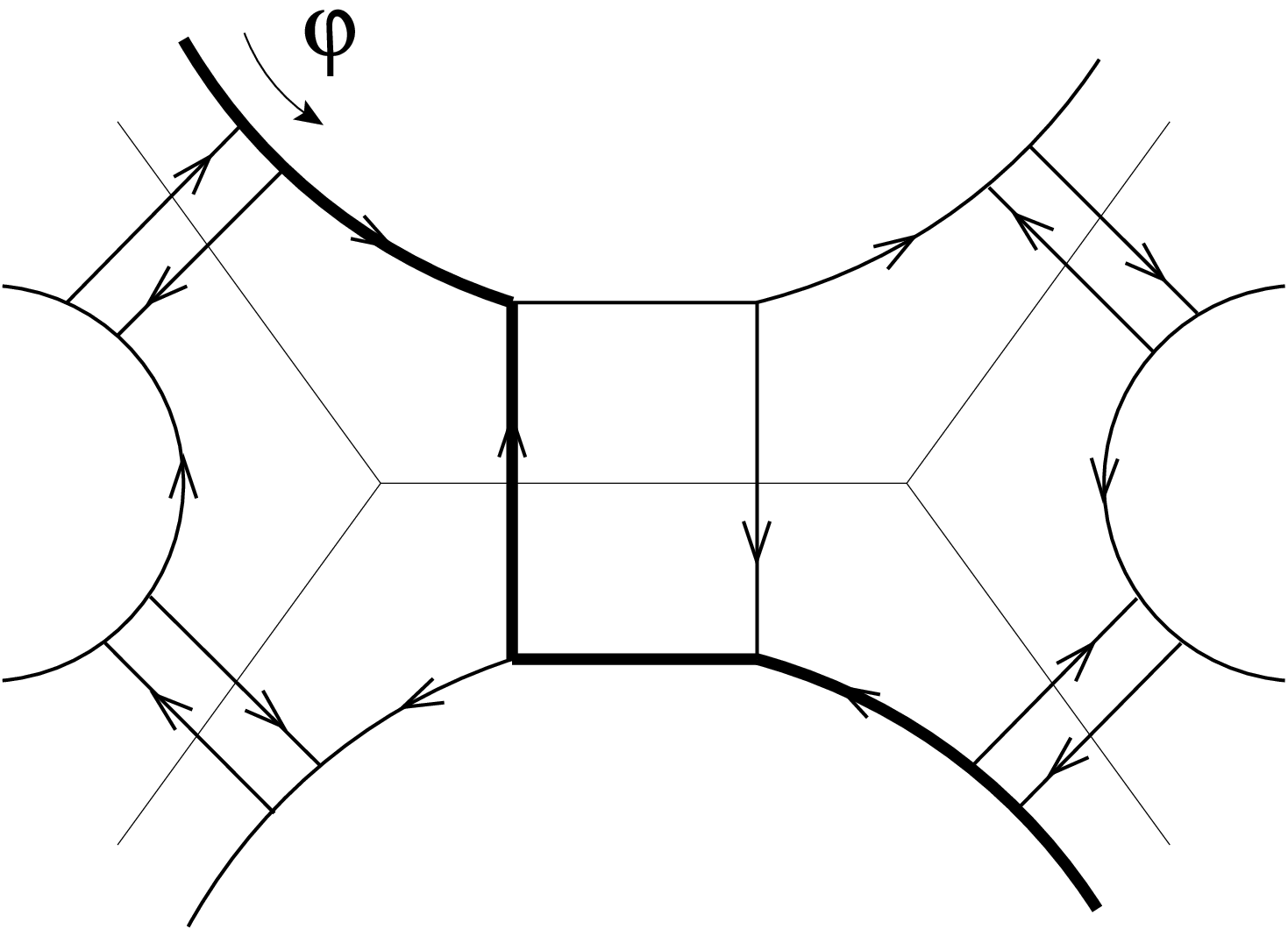}\quad\quad
\includegraphics[width=0.4\textwidth]{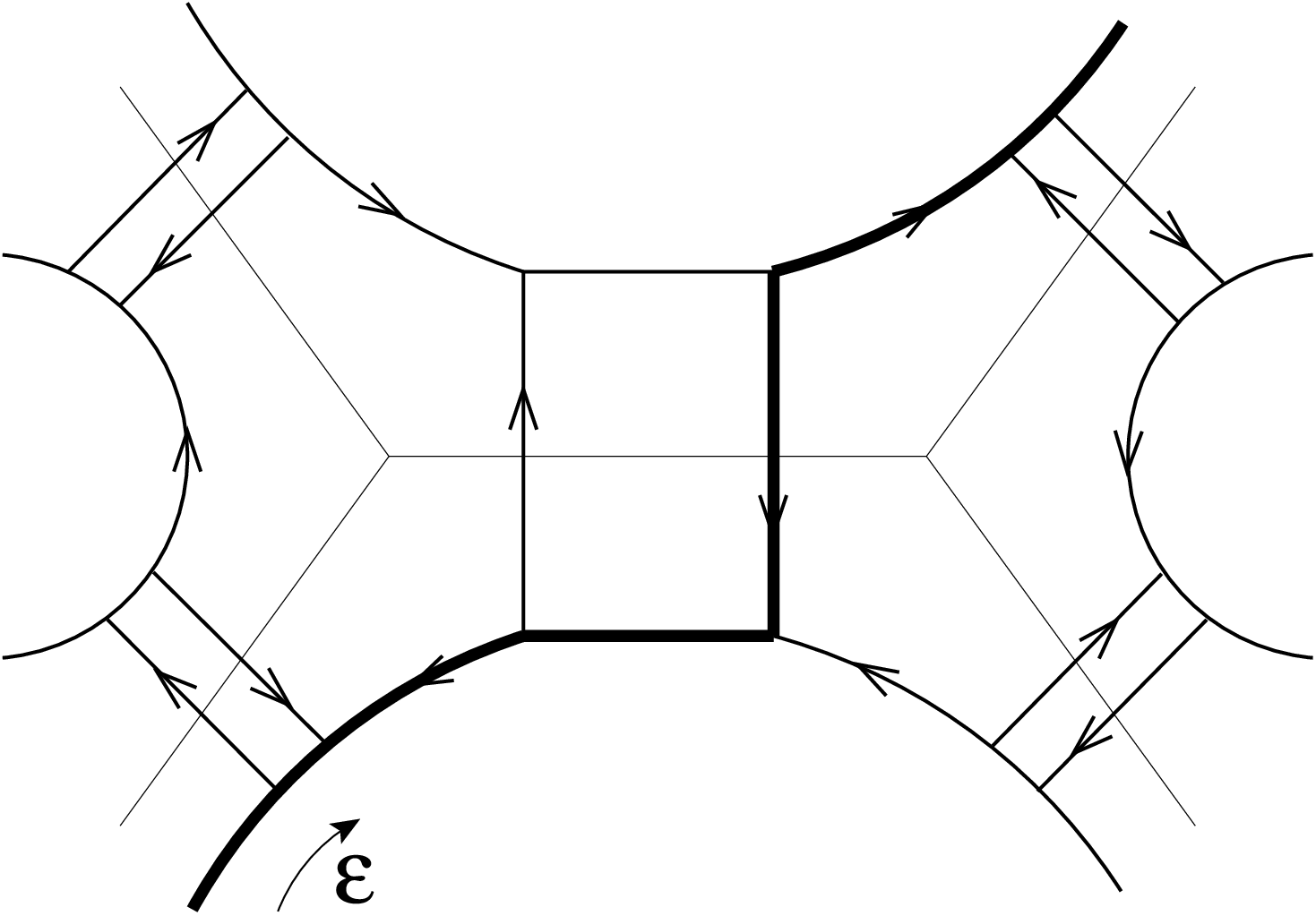}\\         
\includegraphics[width=0.4\textwidth]{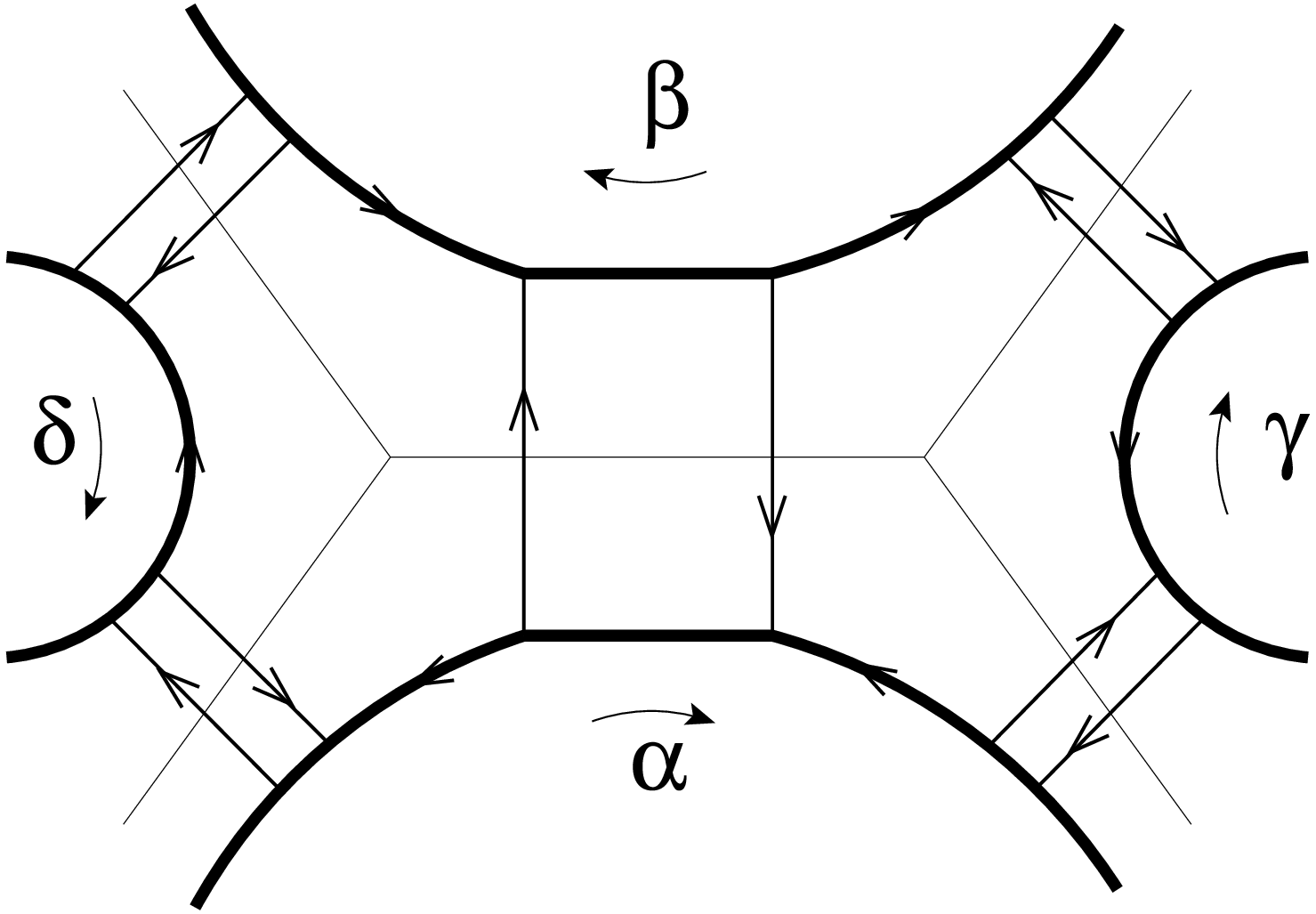}
\caption{The oriented paths $\alpha,\beta,\gamma,\delta,\varepsilon,\phi$ on $\bar{\mathcal{G}}_{\tau_1}$ are indicated by bold lines}  \label{curves}               
\end{figure}

\noindent{\bf Proof.}  Consider a neighborhood of the edge of $\tau_1$ upon which the flip is performed depicted in Figure \ref{curves} where
there are illustrated six distinct oriented paths in $F$ denoted $\alpha, \beta, \gamma, \delta, \epsilon, \phi$, 
which may be completed to closed oriented curves in $F$ and contribute 
to the value of the quadratic form.
A tedious computation given in Appendix II compares formula \rf{q} before and after the flip
for these six paths and determines that there is the unique evolution of orientation class from $\tau_1$ to $\tau_2$ illustrated in Figure \ref{flipper} that leaves invariant these contributions.
\hfill$\blacksquare$

\bigskip

\smallskip

\begin{Thm}\label{oops}
Let $\tau$ be a trivalent fatgraph spine of $F=F^s_g$. Then $\mathcal{O}({\tau})$ is isomorphic to the set $\Omega(F)$ of spin structures on $F$ as affine $H^1(F;\mathbb{Z}_2)$-spaces. 
Moreover,  the action of $MC(F)$ on $\Omega(F)$ lifts naturally to the action of the Ptolemy groupoid on $\mathcal{O}({\tau}$) illustrated in Figure \ref{flipper}.  
\end{Thm}

\smallskip

\noindent By naturality here, we mean that if $\varphi\in MC(F)$, then there is a finite sequence $F\supset\varphi(\tau)=\tau_1\to\tau_2\to\cdots\to\tau_{n}=\tau$ where $\tau_{i+1}$ arises from $\tau_i$ by a flip, for $i=1,\ldots,n-1$.  Suppose that $\omega$ is an orientation on $\tau$ inducing via the identification $\varphi$ the orientation $\omega_1$ on $\tau_1$ and moreover that
the evolution of orientation illustrated in Figure \ref{flipper} serially induces from $\omega_1$
the orientation $\omega'=\omega_{n+1}$ on $\tau_{n+1}=\tau$.  The spin structure on $F$ associated to $\omega$ maps to that of $\omega'$ under the action of $\varphi$ on $\Omega(F)$ by construction.

\section{Coordinates on decorated super-Teichm\"uller space}


Equivalent to the choice of trivalent fatgraph spine $\tau\subset F$ is the specification of
its dual ideal triangulation $\Delta$ of $F$.  An orientation $\omega$ on $\tau$ induces an orientation on $\Delta$ by requiring that the oriented edge of $\Delta$ occurs clockwise from its dual oriented edge in $\tau$ near their point of intersection here using an orientation of the surface $F$.
Dual to the fatgraph reflection at a vertex is the change of orientation on each edge in the frontier of a triangle complementary to $\Delta$.
We consider the lift of $\Delta$ to an ideal triangulation $\tilde\Delta$
of the universal cover $\tilde F\to F$.  Fixing a base point in $F$, the fundamental group $\pi_1=\pi_1(F)$ acts as deck transformations on $\tilde F$ leaving invariant $\tilde\Delta$.
Of course $\tilde F$ is {\sl topologically equivalent} to ${\mathbb D}$, and a hyperbolic metric on $F$ further determines a {\sl metric equivalence}.  In any case, we may consider the collection of ideal vertices $\tilde\Delta_\infty\subset {\mathbb S}^1$ of all the arcs in $\tilde\Delta$ as an abstract set.
\bigskip

\begin{Thm} {\it Fix a surface $F=F_g^s$ of genus $g\geq 0$ with $s\geq 1$ punctures,
where $2g-2+s>0$, and let $\Delta$ be some ideal triangulation of $F$
whose lift $\tilde\Delta$ to the universal cover $\pi:\tilde F\to F$ has ideal vertices
$\tilde \Delta_\infty$.  Suppose that 
$\omega$ is an orientation class on the arcs in $\Delta$ corresponding to a specified spin structure, and assign to each edge of $\Delta$ an even coordinate and to each 
 triangle complementary to $\Delta$ an odd coordinate where the latter are taken modulo an overall sign.  Then there is a function 
$$\ell:\tilde\Delta_\infty\to\hat{L}^+_0$$
uniquely determined up to post-composition with an element of $OSp(1|2)$
so that
\bigskip

\leftskip .3in

\noindent {\rm i)}~if $a,b\in \tilde\Delta_\infty$ span an edge in $\tilde\Delta$, then 
the coordinate on this edge is given by the $\lambda$-length $\sqrt{<\ell(a),\ell(b) >}$;

\medskip

\noindent {\rm ii)}~if $c,b,a\in\tilde\Delta_\infty$ span a  triangle complementary
to $\tilde\Delta$ and occur in the positive order in ${\mathbb S}^1$,
then up to a sign the coordinate of this triangle is given by the
$\mu$-invariant of the positive triple $\ell(c)\ell(b)\ell(a)$; moreover,
if $d,c,a\in\tilde\Delta_\infty$ likewise occur in the positive order
and span a triangle, then the coordinates of these respective triangles
are given up to an overall sign by the $\mu$-invariants of the positive
triples $\ell(c)\ell(b)\ell(a)$ and $\ell(d)\ell(c)\ell(a)$ as related by Proposition \ref{bascalc}.

\leftskip=0ex
\bigskip

\noindent Furthermore, there is a representation $\hat\rho:\pi_1=\pi_1(F)\to OSp(1|2)$ with respect to which $\ell$ is $\pi_1$-equivariant in the sense that $\hat\rho(\gamma) (\ell(a))=\ell(\gamma (a))$ for each $\gamma\in\pi_1$ and $a\in\tilde\Delta_\infty$ so that  
$\pi_1{{{~\atop {\hat\rho}}\atop\to}\atop ~} OSp(1|2){{{~\atop {}}\atop\to}\atop ~} SL(2,{\mathbb R}){{{~\atop {}}\atop\to}\atop ~} PSL(2,{\mathbb R})$ is a Fuchsian representation
whose lift $\pi_1{{{~\atop {\hat\rho}}\atop\to}\atop ~} OSp(1|2){{{~\atop {}}\atop\to}\atop ~} SL(2,{\mathbb R})$ agrees with the specified spin structure.  Moreover $\hat\rho$ is uniquely determined up to conjugacy by an element of $OSp(1|2)$.
}\end{Thm}

\begin{figure}[hbt] 
\centering           
\includegraphics[width=0.75\textwidth]{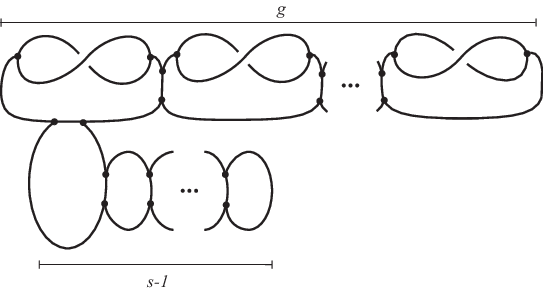}
\caption{A bipartite fatgraph spine in $F_g^s$ for any $g\geq 0$ and $s\geq 1$}  \label{hip}               
\end{figure}

\noindent{\bf Proof.} The argument closely follows the bosonic case \cite{penner,pb} using the putative coordinates to recursively construct the mapping $\ell$ and finally the representation $\hat\rho$.  Let us first consider the case that the trivalent fatgraph spine $\tau\subset F$ dual to $\Delta$ is bipartite, namely, there is a two-coloring of the vertices of $\tau$ so that the endpoints of any edge of $\tau$ have different colors.  An example of a bipartite fatgraph spine  is illustrated in Figure \ref{hip} for each topological type of punctured surface.

Dually letting ${\mathcal T}$ denote the collection of triangles complementary to $\Delta$ in $F$, there is a function $\delta:{\mathcal T}\to\{{\pm 1}\}$ so that for any two triangles $t_1,t_2\in {\mathcal T}$ sharing an edge, we must have $\delta(t_1)\cdot\delta(t_2)=-1$.  Notice that in this bipartite case there are exactly two distinct such functions differing by an overall sign.
The coordinates on $\Delta$ lift to coordinates on $\tilde\Delta$ in the natural way
assigning to each edge $e\in\tilde\Delta$ the even coordinate of $\pi(e)$ and to each
complementary triangle $t\in\tilde {\mathcal T}=\{{\rm complementary~triangles~to}~\tilde\Delta\}$
the odd coordinate $\theta_t$ of $\pi(t)$.  The function $\delta$ likewise lifts to the eponymous function $\delta:\tilde {\mathcal T}\to\{\pm 1\}$ satisfying $\delta(t)=\delta(\pi(t))$.

Choose a distinguished triangle-edge pair or equivalently a distinguished oriented edge of $\tilde\Delta$ which determines the triangle $t\in \tilde {\mathcal T}$
lying to its left in the oriented surface $\tilde F$.  According to Lemma \ref{3pts}, there is a positive triple $ABC$ in the special light cone 
realizing
the putative  $\lambda$-lengths and $\mu$-invariant $\theta_t$
so that the point $B$ is opposite to the distinguished oriented edge, and $ABC$ is furthermore uniquely determined up to fermionic reflection $g^r$.  Choosing a particular representative $ABC$ or $g^r ABC$ at this stage gauge fixes the ${\mathbb Z}_2$-action and thus  determines well-defined signs on all of the fermionic coordinates $\theta_{\bar t}$ for $\bar t\in \tilde {\mathcal T}$.  (In the special case when  the fermionic coordinate on $t$ vanishes, we have $g^r ABC=ABC$ and can still gauge fix at this point to determine the signs on all of the non-vanishing fermionic coordinates.)  This lift to a positive triple in $\hat L^+_0$ of the vertices of $t$ constitutes the basis step of our recursive construction of the mapping $\ell$.

In order to lift the triangle $t'\in \tilde {\mathcal T}$ to the right of the distinguished oriented edge
with its specified fermionic coordinate $\theta_{t'}$,
we shall rely on the basic calculation formulated as Proposition \ref{bascalc} to lift the vertex of $t'$ distinct from those of $t$ to a point $D\in \hat L^+_0$,
however, the coordinates employed to determine $D$ are the specified putative $\lambda$-lengths and the {\sl $\delta$-modified $\mu$-invariant} $\delta(t') \cdot \theta_{t'}$ rather than
$\theta_{t'}$ itself. Conjugating by the $ABC$ prime transformation and its square, we can likewise uniquely lift  to $\hat{L}^+_0$ the vertices of the other two triangles in $\tilde {\mathcal T}$ adjacent to $t$. 

Continue recursively in this way moving each triangle in $\tilde {\mathcal T}$ already lifted into standard position with its specified $\mu$-invariant
using $OSp(1|2)$ and then employing the putative $\lambda$-lengths and $\delta$-modified $\mu$-invariants
to lift the next adjacent triangle in $\tilde {\mathcal T}$ using Proposition \ref{bascalc}
in order to finally uniquely lift all of $\tilde \Delta_\infty$ 
and completely define $\ell:\tilde\Delta_\infty\to\hat L^+_0$ in keeping with
requirements i) and ii) in the statement of the theorem.
There are actually only finitely many values of coordinates employed in this construction since there are only finitely many arcs in $\Delta$ and triangles in ${\mathcal T}$.
In particular, the bodies of $\lambda$-lengths are uniformly bounded above and
below.  As in the pure even case \cite{penner,pb}, it follows from this bounded geometry that the bosonic projection $\tilde \Delta_\infty {{{~\atop \ell}\atop\to}\atop ~} \hat L^+_0 {{{~\atop ~}\atop\to}\atop ~} L^+$ induces a mapping from the ideal triangulation of $\tilde F$ onto a tessellation of all of
${\mathbb D}$.

We claim that changing the distinguished oriented edge of $\tilde \Delta$ used to initiate this construction merely modifies the function $\ell$ by post-composition with a group element in $OSp(1|2)$.
To see this, first notice that if we change orientation of the distinguished oriented edge,
then $t'$ will be
in standard position with its $\mu$-invariant $\theta_{t'}$, and the triangle $t$ is then produced using the $\delta$-modified $\mu$-invariant $\delta({t})\cdot \theta_t$. Applying the switch transformation (and possibly overall fermionic reflection for the ${\mathbb Z}_2$ gauge group), based on the fact that 
$\delta(t)\cdot\delta(t')=-1$ together with Proposition \ref{switch}, we recover the image of the previous construction under the switch transformation, which lies in $OSp(1|2)$, with the correct signs for fermionic parameters. 
Moreover, any change in the distinguished oriented edge can be achieved by means of a finite composition of prime transformations and at most one reversal of orientation. Since prime transformations are again group elements of $OSp(1|2)$, we 
conclude that our equivalence class of lifts under the action of $OSp(1|2)$ is independent of the distinguished oriented edge as required.  Uniqueness of the mapping $\ell$ up to post-composition with an element of $OSp(1|2)$ follows by construction from the uniqueness statements in
Lemma \ref{3pts} and Proposition \ref{bascalc} completing the proof of the first part under the bipartite assumption.

From the lift $\tilde \Delta_\infty {{{~\atop \ell}\atop\to}\atop ~} \hat L^+_0$, we shall presently construct the required representation 
$\hat\rho:\pi_1\to\hat\Gamma<OSp(1|2)$.  To this end,
choose a connected fundamental domain $D\subset \tilde F$ for the action of $\pi_1$ comprised of a collection of triangles in $\tilde {\mathcal T}$ which contains the base triangle $t$
specified earlier which began the recursive construction of $\ell$.  The frontier edges of $D$ in $\tilde F$ arise in pairs
$c,c'$ together with an abstract identification $c'=\gamma(c)$ induced by some $\gamma\in\pi_1$, and we let $c_i'=\gamma_i(c_i)$
enumerate the collection of these edge pairings, for $i=1,\ldots ,4g+2s$, where the $\gamma_i$ thus freely generate $\pi_1$.
To determine the image $\hat\rho(\gamma_i)\in OSp(1|2)$ of $\gamma_i\in\pi_1$ in order to define the representation $\hat\rho$, let us further enumerate two triangles for each edge pairing,
namely, the unique triangles $t_i\supset c_i$ and $t_i'\supset c_i'$ in $\tilde {\mathcal T}$
where $t_i\subset D$ and $t_i'\not\subset D$.  Notice as before that the coordinates for $t_i$ or $t_i'$
have gauge-fixed signs from the basis step of the construction of $\ell$, and the parameters on $t_i$ lifted from $F$ agree with those
on $t_i'$ since $\gamma_i(t_i)=t_i'$ implies $\pi(t_i)=\pi(t_i')$.


Let $T_i=\ell(t_i)$ be the corresponding positive triple in $\hat L^+_0$. 
According to Lemma \ref{3pts} again, there is a unique $h_i\in OSp(1|2)$ such that $T_i\cdot h_i=\bar T_i$, where $OSp(1|2)$ acts on the right and
$\bar T_i$ indicates the triple of points in standard position realizing the $\lambda$-lengths and signed $\mu$-invariant of $T_i$
so that the fermion-dependent vertex is opposite to $c_i$. 
Let $h'_i\in OSp(1|2)$ be the analogous group element for $T'_i$.
We must furthermore take account of the spin structure in accordance with 
\cite{natpaper,natanzon} relating the sign of the trace of the bosonic reduction of $\hat\rho(\gamma)$ with the value of the corresponding quadratic form on the underlying mod two homology class $[\gamma]\in H_1(F;{\mathbb Z}_2)$, for $\gamma\in\pi_1$.  
In particular, notice that
the fermionic reflection $g^r$ has the central element $-I\in SL(2,{\mathbb R})$ as its bosonic reduction.
Given the quadratic form $q\in{\mathcal Q}(F)$ corresponding via Theorem \ref {qando} to our specified orientation $\omega$, we define
\begin{eqnarray}\nonumber
\hat\rho(\gamma_i)=\nonumber\begin{cases}
&h_i^{-1}\cdot g^r\cdot h'_i,\\ 
&{\rm if}~~ {\rm trace}(h_i^{-1}\cdot h'_i)>0 ~\&~ q([\gamma_i])=0 ~~ {\rm or}~~  {\rm trace}(h_i^{-1}\cdot h'_i)<0~\&~ q([\gamma_i])=1;\\
&\\
&h_i^{-1}\cdot h'_i, \\
&{\rm if}~~ {\rm trace}(h_i^{-1}\cdot h'_i)<0~\&~ q([\gamma_i])=0 ~~ {\rm or}~~  {\rm trace}(h_i^{-1}\cdot h'_i)>0~\&~ q([\gamma_i])=1.
\end{cases}
\end{eqnarray}

By construction, these group elements compose correctly so as to produce a representation, and since $\pi_1$ is a free group,
there are no relations to confirm whence $\hat\rho$ is indeed a representation of $\pi_1$ in $OSp(1|2)$ onto a subgroup $\hat\Gamma=\hat\rho(\pi_1)< OSp(1|2)$.
Also by construction, the projectivized bosonic reduction $\Gamma=\rho(\pi_1)<PSL(2,{\mathbb R})$
leaves invariant the tessellation of ${\mathbb D}$ discussed before, and an argument 
in \cite{penner,pb} going back to Poincar\'e thus proves that $\Gamma$ is
indeed a Fuchsian group uniformizing the punctured surface.  Finally, the bosonic reduction $\tilde\Gamma=\tilde\rho(\pi_1)<SL(2,{\mathbb R})$ itself as a lift of $\Gamma<PSL(2,{\mathbb R})$ gives the correct spin structure on the underlying Riemann surface in  keeping with \cite{natpaper,natanzon} by construction since multiplication by $g^r$ alters the sign of the trace of the bosonic reduction.

This completes the required construction of a representation from the asserted parameters including the spin structure. The procedure described is clearly equivariant for the fermionic reflection under the initial choice of sign on the $\mu$-invariant in the base triangle, and
the distinction between pairs $t_i\subset D$ and $t_i'\not\subset D$ amounts only to replacing a generator by its inverse. 
There are thus two essential
choices: the choice of base triangle $t$ to begin the inductive construction as well as the fundamental domain $D$ containing it.  
We must show that these choices are resolved by the overall conjugacy in the definition of the decorated super-Teichm\"uller space.

To this end, consider two such fundamental domains containing specified base triangles $t_i\subset D_i$ with corresponding mappings $\ell_i:\tilde\Delta_\infty\to\hat L^+_0$
and representations $\hat\rho_i:\pi_1\to OSp(1|2)$, for $i=1,2$.
There is a unique triangle $ t_1'\subset D_2$ whose projection to $F$ agrees with that of $t_1$ since $D_2$ is a fundamental domain for the action of $\pi_1$, whence
$t_1$ and $\bar t'_1$ thus share the same invariants.  By Lemma \ref{3pts} again ${\mathbb Z}_2$ gauge-fixed as before, there is a unique
$g\in OSp(1|2)$ so that $g(t_1)=t_1'$.  By uniqueness in Lemma \ref{3pts} and Proposition \ref{bascalc}, we must then have $\ell_1=\ell_2\circ g$, and so
$g$ furthermore conjugates $\hat\rho_1$ to $\hat\rho_2$ as required.  
Notice that the trace of the bosonic reduction of an element of $OSp(1|2)$ is invariant under conjugacy, so the spin structure is left invariant by conjugation of representations.

We must finally extend the construction from bipartite to general trivalent fatgraph spines $\tau\subset F$.  To this end according Corollary D, there is a finite sequence of flips
starting from $\tau$ and ending with a bipartite fatgraph spine $\tau'\subset F$ such as the one depicted in Figure \ref{hip}
to which we may apply the construction just described based upon the coordinates
on $\tau'$ induced from those on $\tau$ via super Ptolemy transformations in Proposition \ref{pt*} and Theorem \ref{fl1}  as well as the orientation class on $\tau'$ induced from that on $\tau$ in Theorem \ref{oops} determining the spin structure.  Since these super Ptolemy transformations are computed relative to a fixed configuration of points in $\hat L^+_0$, this gives a well-defined lift $\ell:\tilde\Delta_\infty\to\hat L^+_0$ for any trivalent fatgraph spine and hence a corresponding representation $\hat\rho:\pi_1\to OSp(1|2)$ which is determined up to conjugacy since the mapping $\ell$ is determined up to post-composition with an element of $OSp(1|2)$.  
\hfill $\blacksquare$\\

\begin{Cor}
Fix a trivalent fatgraph spine $\tau\subset F$ for a surface $F=F_g^s$ of negative Euler characteristic and specify an orientation on the edges of $\tau$ determining the 
component $C$ of $S\tilde T(F)$.  Then $\lambda$-lengths on the edges  together with $\mu$-invariants on the vertices of $\tau$, the latter taken modulo an overall sign, provide global coordinates on $C$.  Moreover, these coordinates are natural in the sense
that if $\varphi\in MC(F)$ and $\varphi(\tau)=\tau_1-\tau_2-\cdots -\tau_n=\tau$ is a sequence of trivalent fatgraph spines of $F$ with consecutive ones related by a flip, then we identically induce coordinates
and orientation class on $\varphi(\tau)$ from these data on $\tau$ using $\varphi$ and perform the corresponding sequence of super Ptolemy transformations in Proposition \ref{pt*} and Theorem \ref{fl1}
and evolution of orientation in Theorem \ref{oops}
to induce new coordinates and orientation on $\tau$ itself.  Then these induced coordinates and orientation class on $\tau$ describe the action of $\varphi$ on $S\tilde T(F)$.
\end{Cor}

\medskip

\noindent{\bf Proof.}
This follows directly from the previous result and Theorem \ref{oops} together with the observation that $\lambda$-lengths and $\mu$-invariants are defined intrinsically in ${\mathbb R}^{2,1|2}$.
\hfill$\blacksquare$

\section{Shear coordinates and Ptolemy-invariant 2-form}

First, let us recall in some detail the Ptolemy transformations studied in Section \ref{sec:bcpt}. There were two aspects to the calculation corresponding to even and odd,
namely, the even Ptolemy transformation is a simple modification
\begin{eqnarray}\label{ptle}
 ef=(ac+bd)\Big(1+\frac{\sigma\theta\sqrt{Z}}{1+Z}\Big)
\end{eqnarray}
 of the standard pure even case $ef=ac+bd$, and
the odd Ptolemy transformation again in the notation of Section \ref{sec:bcpt} is given by
\begin{eqnarray}\label{0opt}
\nu=\frac{\sigma+\theta\sqrt{Z}}{\sqrt{1+Z}},\quad
\mu=\frac{\sigma\sqrt{Z}-\theta}{\sqrt{1+Z}}.
\end{eqnarray}
Taken together \rf{ptle} and \rf{0opt} give expression to the {\it super Ptolemy transformation}.

These formulas
in particular
specialize to those of \cite{zababaher}
where shear coordinates on the super-Teichm\"uller space are introduced.
Namely as in the classical case, given an ideal triangulation $\Delta$ of the surface $F$, to each edge $AC$ as in Figure \ref{figint} is associated its cross ratio $Z_e=Z=\frac{ac}{bd}$, written again in terms of $\lambda$-lengths.  Given the two positive triples $CBA$ and $DCA$ in the special light cone
$\hat{L}^+_0$, denote by $Z_a$, $Z_b$, $Z_c$, $Z_d$ these cross ratios associated to each respective edge $a,b,c,d$ also illustrated in Figure \ref{figint}
identifying an edge with its $\lambda$-length here for convenience; of course, each of these depends in turn on the $\lambda$-lengths of other nearby edges, for instance,
$Z_b=\frac{eg}{ah}$, where $g$, $h$ are $\lambda$-lengths on the frontier edges of the other complementary triangle to $\Delta$ than $CBA$ that contains $CB$.  Nevertheless, under the even Ptolemy transformation, we find that $Z_b=\frac{eg}{ah}$ transforms to
\begin{eqnarray}
&&Z_b \frac{ac}{ef}=Z_b(1+Z^{-1})^{-1}\Big(1+\sigma\theta(Z^{\frac{1}{2}}+Z^{-\frac{1}{2}})^{-1}\Big)^{-1}=
Z_b(1+Z^{-1}+\sigma\theta\sqrt{Z^{-1}})^{-1},\nonumber
\end{eqnarray}
and one similarly computes
\begin{eqnarray}
&&Z_a\mapsto Z_a(1+Z+\sigma\theta\sqrt{Z}), \quad Z_c\mapsto Z_c(1+Z+\sigma\theta\sqrt{Z}),\nonumber\\
&&Z_d\mapsto Z_d(1+Z^{-1}+\sigma\theta\sqrt{Z^{-1}})^{-1}.
\end{eqnarray}
In fact, these transformations together with equation \rf{0opt} coincide with those in \cite{zababaher} up to a conventional inversion of the cross ratio coordinate (see Figure 2 in \cite{zababaher}).

In the pure even case  \cite{penner-jdg2,pb}, the K\"ahler 2-form of the
Weil-Petersson Hermitian metric
on the Teichm\"uller space $T(F)$ or moduli space $M(F)$ of a punctured surface $F$
was computed relative to the $\lambda$-lengths on any convenient triangulation
$\Delta$ of $F$.  
Namely, it is given by
\begin{eqnarray}\label{2fbos}
\omega_{\Delta}=2~\sum d\log a\wedge d\log b+d\log b\wedge d\log c
+d\log c\wedge d\log a,
\end{eqnarray}
where the sum is over all complementary triangles $T$ to $\Delta$ in $F$
with frontier edges of $T$ occurring in the cyclic order $a,b,c$ compatibly with 
the clockwise orientation on $T\subset F$.  It is not difficult to check directly
that this expression is invariant under the pure even Ptolemy transformation $ef=ac+bd$, and indeed in the 1980's before it was recognized as this particular 2-form, it was nevertheless already confirmed to be invariant under Ptolemy transformations and hence arise from some 2-form on the quotient moduli space $T(F)/MC(F)$.  Much this same computation applies in the current case.

\begin{Thm}  If $\Delta$ is an ideal triangulation of the punctured surface $F$,
then consider the even 2-form on
$S\widetilde{ \mathcal T}(F)$ given by
\begin{eqnarray}\label{2fsuper}
\hat\omega_\Delta=
\sum d\log a\wedge d\log b+d\log b\wedge d\log c
+d\log c\wedge d\log a -(d\theta)^2,\nonumber
\end{eqnarray}
where the sum is over all triangles whose consecutive edges in the clockwise ordering have $\lambda$-lengths $a, b, c$ and $\mu$-invariant $\theta$.
Then this 2-form is invariant under super Ptolemy transformations.
\end{Thm}

\medskip

\noindent Notice that  $\hat\omega_\Delta$ is manifestly invariant under the fermionic reflection changing the signs of all fermions, and we have dropped the pre factor 2 here compared to the K\"ahler form so that if one converts the expression $\hat\omega_\Delta$ from $\lambda$-lengths to shear coordinates as described before, then the resulting 2-form is associated with the Poisson bracket given in \cite{zababaher}.

\medskip

\noindent{\bf Proof.} Let us remind \cite{berezin}  the reader that on supermanifolds the de Rham operator anticommutes with odd constants, and in local coordiantes $(\{x_i\}, \{\theta_j\})$  is given by $d=dx_i\p_{x_i}+d\theta_i\p_{\theta_i}$, where the odd derivative acts from the left and the $d\theta_i$ are even. 

Adopt the notation $\tilde{x}=d\log x=\frac{dx}{x}$ for any invertible expression $x$ and compute as in the pure even case \cite{penner-jdg2,pb}
for a pair $CBA$ and $DCB$ of positive triples that the contribution to $\hat\omega_\Delta$ before the flip
on $AC$ in the notation of Figure \ref{flip0}
is given by
\begin{eqnarray}\label{qbo1}
\t a\t b+\t b\t e+\t e\t a+\t e\t d+\t d\t c+\t c\t a-(d\theta)^2-(d\sigma)^2
\end{eqnarray}
and after the flip is given by
\begin{eqnarray}\label{qbo2}
\t b\t c+\t c\t f+\t f \t b+\t f \t d+\t d\t a +\t a\t f-(d\nu)^2-(d\mu)^2.
\end{eqnarray}   
We must show that these two expressions coincide.

To this end, notice that if $U=\frac{\sigma\theta\sqrt{Z}}{1+Z}=
\frac{\sigma\theta}{Z^{\frac{1}{2}}+Z^{-\frac{1}{2}}}$, then 
\begin{eqnarray}
\t e+\t f=\frac{1}{ac+bd} (ac(\t a+\t c)+bd(\t b+\t d)) +d U,\nonumber\\
\end{eqnarray} 
and we find
$$\begin{aligned}
\t b\t c+\t c\t f+\t f \t b+\t f \t d+\t d\t a +\t a\t f&=\t b\t c+\t d \t a+\t f(\t b+\t d-\t a -\t c)\\
&=\t b\t c+\t d \t a-\t e(\t b+\t d-\t a -\t c)\\
&~~~~+\frac{ac}{ac+bd}(\t a+\t c)(\t b+\t d)+\frac{bd}{ac+bd}(\t a+\t c)(\t b+\t  d)- dU\t Z,
\end{aligned}$$
where
\begin{eqnarray}\label{du}
dU\t Z=\frac{d(\sigma\theta)dZ}{(1+Z)\sqrt{Z}}.
\end{eqnarray}
Meanwhile, the super Ptolemy transformation gives
  $$\begin{aligned}
(d\mu)^2+(d\nu)^2&=d\sigma^2+d\theta^2+2\frac{d\sigma}{\sqrt{1+Z}}d\Big(\frac{\sqrt{Z}}{\sqrt{1+Z}}\Big)\theta+
2\frac{\sqrt{Z}d\theta}{\sqrt{1+Z}}d\Big(\frac{1}{\sqrt{1+Z}}\Big)\sigma\\
&=-2\frac{d\sigma\sqrt{Z}}{\sqrt{1+Z}}d\Big(\frac{1}{\sqrt{1+Z}}\Big)\theta-2\frac{d\theta}{\sqrt{1+Z}}d\Big(\frac{\sqrt{Z}}{\sqrt{1+Z}}\Big)\sigma\\
&=d\sigma^2+d\theta^2+\frac{d(\theta\sigma)dZ}{(1+Z)\sqrt{Z}}\\
\end{aligned}$$
which coincides with equation \rf{du} as required.\hfill $\blacksquare$

\section*{Appendix I. $OSp(1|2)$: Notation and conventions}
In this appendix, we provide basic information concerning the Lie supergroup $OSp(1|2)$ and its Lie superalgebra. 
Our cursory treatment here is presumably sufficient for the purposes of the text, and we refer the interested reader to \cite{berezin,manin}
for further details about general Lie superalgebras and supergroups.

Let us first introduce certain conventions (which differ from those in \cite{zababaher}).
Given a Lie superalgebra $\mathfrak{g}$, one can consider its Grassmann envelope, namely, 
the Lie superagebra $\mathfrak{g}(S)=S\otimes\mathfrak{g}$ for some Grassmann algebra $S$ with decomposition $S=S_0\oplus S_1$
into even and odd elements.
It follows that $\mathfrak{g}(S)$ is both a right and left $S$-module, i.e.,
$s\otimes T=(-1)^{|s||T|}(1\otimes T)(s\otimes 1)$ if $s\in S$ and $T\in\mathfrak{g}$ are homogeneous elements of respective degrees $|s|$ and $|T|$.
This rule allows one to construct a representation of the corresponding Lie superalgebra $\mathfrak{g}(S)$ in the space $S\otimes\mathbb{R}^{m|n}$ from a given representation of $\mathfrak{g}$ in $\mathbb{R}^{m|n}$. One can then produce a representation of the corresponding Lie group $G(S)$ by exponentiating pure even elements from $\mathfrak{g}(S)$ in $S\otimes\mathbb{R}^{m|n}$.

When writing a super matrix
$\left(\begin{array}{cc}
A & B  \\
C & D \\
\end{array} \right)$
representing the action of $G(S)$ or $\mathfrak{g}(S)$ as elements of $S\otimes End(\mathbb{R}^{m|n})$ on $S\times\mathbb{R}^{m|n}$, 
notice that for pure even supermatrices, the composition rule is given by
\begin{eqnarray}
\left(\begin{array}{cc}
A & B  \\
C & D \\
\end{array} \right)
\left(\begin{array}{cc}
A' & B'  \\
C' & D' \\
\end{array} \right)=\left(\begin{array}{cc}
AA'-BC' & AB'+BD'  \\
CA'+DC' & DD'-CB' \\
\end{array} \right),
\end{eqnarray}
where  the products on the right hand side are the usual products of (super)matrices. 

The usual (super)matrix multiplication (without the minus signs above) is recovered upon replacing $B$ with $-B$. This difference in sign is related to the fact that one typically considers the action of group elements on $\mathbb{R}^{m|n}_S={S_0}^{\times m}\times {S_1}^{\times n}$, which can be identified with the space of even elements in $S\otimes \mathbb{R}^{m|n}$, and the extra minus sign in front of $B$ comes from that isomorphism.
However, throughout this paper we keep the above convention for multiplication of superalgebras (with the extra signs) since it gives a cleaner
relationship with the representation of the original Lie superalgebra .

Another essential ingredient is the superdeterminant or Berezinian of an even supermatrix $M= \bigl(\begin{smallmatrix}
A&B\\ C&D\\
\end{smallmatrix} \bigr)$ with $D$ invertible which is defined as
\begin{eqnarray}
sdet(M)=(det{D})^{-1}(A+BD^{-1}C),
\end{eqnarray}
where the products of matrices are standard matrix products,
and we have the nonstandard plus sign $A+BD^{-1}C$ in the definition
reflecting our conventions.

Now, the Lie superalgebra $OSp(1|2)$ has three even $h, X_{\pm}$ and two odd generators $v_{\pm}$ satisfying the commutation relations
\begin{eqnarray}
[h,v_{\pm}]=\pm v_{\pm}, \quad[v_{\pm},v_{\pm}]=\mp 2X_{\pm},\quad  
[v_+,v_{-}]=h.
\end{eqnarray}
The corresponding realization via $(2|1)\times (2|1)$ supermatrices is given by
\begin{eqnarray}
v_{+}=\left( \begin{array}{ccc}
0 & \hskip1ex0 & 1 \\
0 & \hskip1ex0 & 0 \\
0 & -1 & 0 \end{array} \right), \quad 
v_{-}=\left( \begin{array}{ccc}
0 & 0 & 0 \\
0 & 0 & 1 \\
1 & 0 & 0 \end{array} \right), \quad 
h=\left( \begin{array}{ccc}
1 & \hskip1ex0 & 0 \\
0 & -1 & 0 \\
0 & \hskip1ex0 & 0 \end{array} \right) ,\nonumber
\end{eqnarray}
and the corresponding supergroup $OSp(1|2)$ can be faithfully realized as $(2|1)\times (2|1)$ supermatrices $g$ with $sdet$ equal to one obeying the relation
\begin{eqnarray}
g^{st}Jg=J,
\end{eqnarray}
where 
\begin{eqnarray}
J=\left( \begin{array}{ccc}
~0 & 1 & ~0 \\
-1 & 0 & ~0 \\
~0 & 0 & -1 \end{array} \right)
\end{eqnarray}
and where the supertranspose $g^{st}$ of $g$ is given by
\begin{eqnarray}
g=\left( \begin{array}{ccc}
a & b & \alpha \\
c & d & \beta \\
\gamma & \delta & f \end{array} \right)\quad
{\rm implies} \quad
g^{st}=\left( \begin{array}{ccc}
~a & ~c & \gamma \\
~b & ~d & \delta \\
-\alpha & -\beta & f \end{array} \right).
\end{eqnarray}
This leads to the following system of constraints on the entries of $g$:
\begin{eqnarray}
&&f=1+\alpha\beta, \quad ad-bc=f^{-1}, \quad d\gamma-c\delta=\beta\nonumber\\
&&\delta=b\beta-d\alpha, \quad \gamma=a\beta-c\alpha, \quad \alpha=b\gamma-a\delta.
\end{eqnarray}
In particular, $(\begin{smallmatrix}a&b\\c&d\\\end{smallmatrix})\mapsto
( \begin{smallmatrix}a&b&0\\c&d&0\\0&0&1\\\end{smallmatrix})$ describes the canonical inclusion $SL(2,{\mathbb R})<OSp(1|2)$, and  the body of the upper-left
$2\times 2$ block conversely gives the bosonic reduction in $SL(2,{\mathbb R})$ of a matrix in $OSp(1|2)$.

\section*{Appendix II. Flip action on oriented fatgraphs}

\begin{figure}[hbt] 
\centering           
\includegraphics[width=0.5\textwidth]{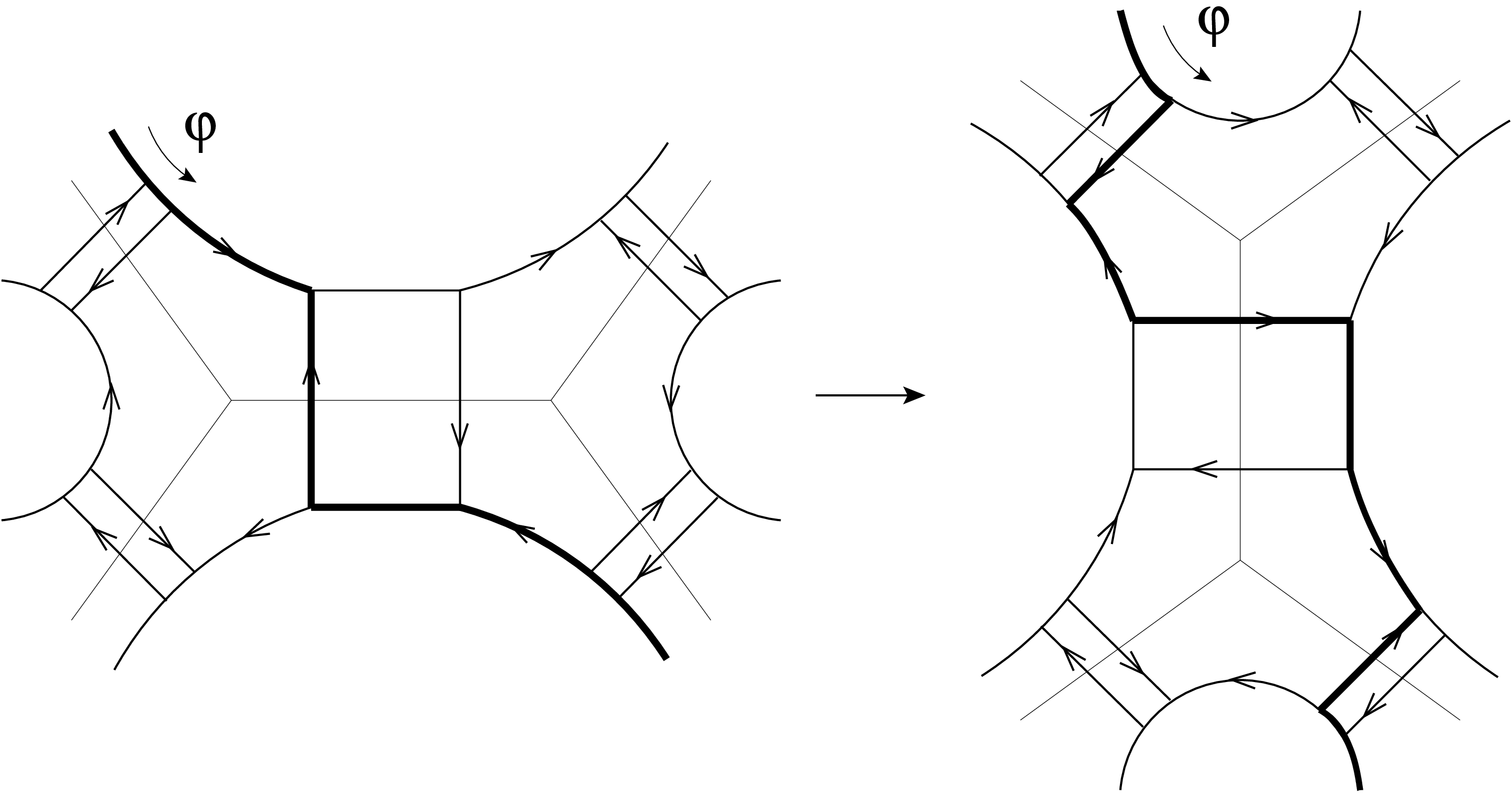}\\
\vspace*{4mm}
\includegraphics[width=0.5\textwidth]{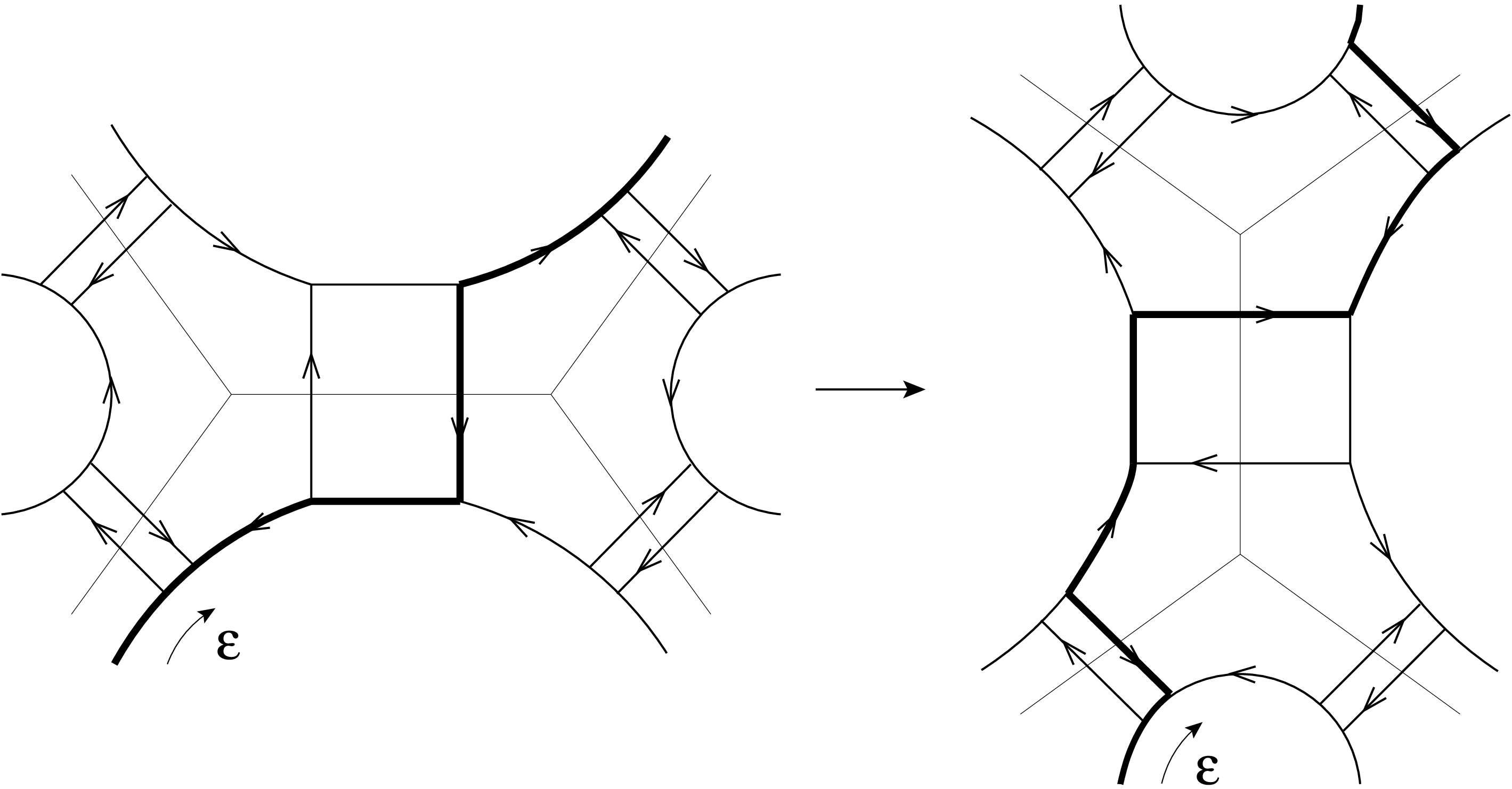}
\caption{Arcs $\varphi$ and $\varepsilon$ before and after the flip}  \label{curveflip}               
\end{figure}

In this appendix, we prove the last part of Theorem \ref{oops} and compute the evolution of orientations under flips illustrated in Figure \ref{flipper}. We must verify that the contributions of the arcs $\alpha,\beta,\gamma, \delta, \varepsilon,\varphi$ illustrated in
Figure \ref{curves} to a quadratic form via Theorem \ref{cr} with the canonical dimer ${\mathcal D}$ for the specified orientations of edges is the same before and after the flip. 
To that end for the purpose of explicit computation, let us fix the representative arcs $\varphi$, $\varepsilon$ as edge-paths before and after the flip as illustrated in Figure \ref{curveflip}.
Representatives edge-paths of the other arcs $\alpha,\beta,\gamma, \delta $ are taken to lie entirely in the boundary.

Now in the notation of Theorem \ref{cr} for the canonical dimer ${\mathcal D}$, we find
$\ell^{\mathcal{D}}_{\alpha}=\ell^{\mathcal{D}}_{\beta}=
\ell^{\mathcal{D}}_{\gamma}=
\ell^{\mathcal{D}}_{\delta}=0 ~(\rm{mod}~2)$ both before and after the flip,  $\ell^{\mathcal{D}}_{\varphi}=\ell^{\mathcal{D}}_{\varepsilon}=1 ~(\rm{mod}~2)$ before the flip and $\ell^{\mathcal{D}}_{\varphi}=1 ~(\rm{mod}~2)$, $\ell^{\mathcal{D}}_{\varepsilon}=0 ~(\rm{mod} ~2)$ after the flip. The numbers $n_C^K$ of Theorem \ref{cr} are tabulated in Figure~10 in the eight cases where the interior edge runs from left to right
with the northwest leaf pointing towards it as may always be arranged by fatgraph reflections on the two vertices.

Uniqueness of this solution for the orientations after the flip is obvious since changing
orientation on anything other than all leaves simultaneously evidently
changes certain of the numbers $n_C^K$.  Moreover, changing the orientation on
each leaf is a composition of the fatgraph reflections at the two interior vertices.

Finally, enumerating the eight cases in Figure~10 in the manner $\begin{smallmatrix}1&2\\3&4\\5&6\\7&8\\\end{smallmatrix}$, one finds by inspection that cases 1, 3, 4, 8 are identical as are cases 2, 5, 6, 7 thus leading to the two cases illustrated in Figure \ref{flipper}.\hfill$\blacksquare$

\begin{figure}[hbt]\label{flipall}
\centering
\includegraphics[width=0.6\textwidth]{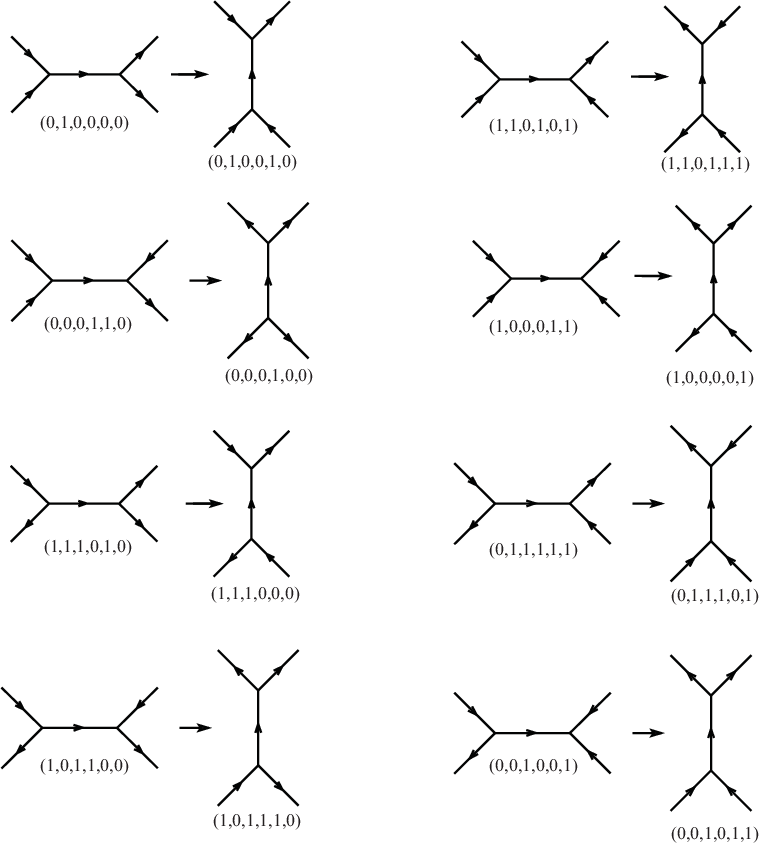}
\caption{All possible salient cases of orientation before the flip
together with resulting orientation after. 
The numbers $n_C^K$ of Theorem \ref{cr} are
tabulated as a vector $(n_\alpha^K,n_\beta^K,n_\gamma^K,n_\delta^K,n_\varepsilon^K,n_\varphi^K)$ in each case.}
\end{figure}


\begin{thebibliography}{99}
\bibitem{beardon} A.\ Beardon, The Geometry of Discrete Groups, Springer, 1983, MR0474012, Zbl 0528.30001.
\bibitem{berezin}F.\ A.\ Berezin, {\it Introduction to superanalysis}, Springer, 1987, MR914369, Zbl 0659.58001.
\bibitem{zababaher} F.\ Bouschbacher, {\it Shear coordinates on the super-Teichm\"uller space}, PhD thesis, Strasbourg, 2014, https://tel.archives-ouvertes.fr/tel-00835500v1/document
\bibitem{BCF} F.\ Bouschbacher, V.\ Fock, F.\ Constantino,
{\it Shear coordinates on the super-Teichm\"uller space}, to appear.
\bibitem{bh}P.\ Bryant, L.\ Hodgkin, {\it Nielsen's theorem and super-Teichm\"ullerspace}, Ann.\ Inst.\ H.\ Poincare Phys.\ Theor., {\bf 58} (3) (1993) 247-265, MR1222942,  Zbl 0787.32013.  
\bibitem{Fock-Chekhov} L.\ Chekhov, V.\ V.\ Fock, {\it A quantum Teichm\"uller space}, 
Theor.Math.Phys. {\bf 120} (1999) 1245-1259; Teor.\ Mat.\ Fiz.\ {\bf 120} (1999) 511-528,  arXiv:math/9908165, MR1737362, Zbl 0986.32007.
\bibitem{cr1}D. Cimasoni, N. Reshetikhin, {\it Dimers on surface graphs and spin structures.\ I}, Comm.\ Math.\ Phys.\ {\bf 275} (2007), 187-208, MR2335773, Zbl 1135.82006.
\bibitem{cr2}D. Cimasoni, N. Reshetikhin, {\it Dimers on surface graphs and spin structures.\ II}, Comm.\ Math.\ Phys.\ {\bf 281} (2008), 445-468, MR2410902, Zbl 1168.82012.
\bibitem{crane-rabin} L.\ Crane, J.\ Rabin, {\it Super Riemann surfaces: uniformization and Teichm\"uller theory}, Comm.\ Math.\ Phys.\ {\bf 113} (1988) 601-623, MR0923633, Zbl 0659.30039.
\bibitem{schwarz} S.N. Dolgikh, A.A. Rosly, A.S. Schwarz, {\it Supermoduli spaces}, Commun. Math. Phys. 135 (1990) 91-100, MR1086753, Zbl 0715.32008.
\bibitem{DW}R. Donagi, E. Witten, {\it Supermoduli Space Is Not Projected}, String-Math 2012, 19–71, Proc. Sympos. Pure Math., 90, Amer. Math. Soc., Providence, RI, 2015, arXiv:1304.7798, MR3409787, Zbl 1356.14021.
\bibitem{flp} A.\ Fathi, F.\ Laudenbach, V.\ Poenaru,
Travaux de Thurston sur les surfaces, {\it Asterisque} (1979){\bf 66-67}, MR0568308,  Zbl 0731.57001. 
\bibitem{ford} L.\ Ford, Automorphic Functions, American Mathematical Society (1951).
\bibitem{hodgkin}L.\ Hodgkin,{\it Super Teichm\"uller Spaces: Punctures and Elliptic Points}, 
Lett.\ Math.\ Phys.\ {\bf 15} (1988) 159-163, MR0943988.  
\bibitem{johnson} D.\ Johnson, {\it Spin structures and quadratic forms on surfaces}, J.\ London Math.\ Soc.\ {\bf 22} (1980) 365-373, MR0588283, Zbl 0454.57011. 
\bibitem{Kashaev} R. Kashaev, {\it Quantization of Teichm\"uller spaces and the quantum dilogarithm}, Lett.\ Math.\ Phys.\ {\bf 43} (1998) 105-115, MR1607296, Zbl 0897.57014.
\bibitem{manin} Yu.\ I.\ Manin, {\it Topics in noncommutative geometry},  Princeton University Press, 1991, MR1095783, Zbl 0724.17007.
\bibitem{milnor} J.\ Milnor, Remarks concerning spin manifolds in {\it Differential and Combinatorial Topology}, Princeton, 1965, MR0180978, Zbl 0132.19602.
\bibitem{natpaper}S.M. Natanzon, {\it Moduli of Riemann surfaces, Hurwitz-type spaces, and their superanalogues}, Russian Math.\ Surveys {\bf 54} 1 (1999) 61-117, MR1706839, Zbl 0948.32018.
\bibitem{natanzon} S.M. Natanzon, {\it Moduli of Riemann Surfaces, Real Algebraic Curves, and Their Superanalogs}, Transl.\ Math.\ Monographs {\bf 225},
American Mathematical Society, 2004, MR2075914, Zbl 1056.14033.
\bibitem{penner}R.\ C.\  Penner, {\it The decorated Teichm\"uller space of punctured surfaces}, 
Comm.\ Math.\ Phys.\ {\bf 113} (1987), 299-339,   Zbl 0642.32012, MR0919235.
\bibitem{penner-jdg2} R.\ C.\  Penner, {\it Weil-Petersson volumes}, J.\ Diff.\ Geom.\ {\bf 35} (1992), 559-608, MR1163449, Zbl 0768.32016.
\bibitem{pb} R.\ C.\  Penner, {\it Decorated Teichm\"uller theory}, European Mathematical Society, 2012, MR3052157, Zbl 1243.30003. 
\bibitem{strebel} K. Strebel, {\it Quadratic Differentials}, Ergebnisse der Mathematik und ihrer Grenzgebiete {\bf 5} Springer-Verlag, Berlin, 1984, MR0743423, Zbl 0547.30001.
\bibitem{verlinde} H. Verlinde,{\it Conformal field theory, two-dimensional quantum gravity and quantization of Teichm\"uller space}, Nucl.\ Phys.\ B {\bf 337} (1990) 652-680, MR1057726.
\bibitem{witten} E. Witten, {\it Notes On Super Riemann Surfaces And Their Moduli},  arXiv:1209.2459. 


\end{thebibliography}
\end{document}